\documentstyle[amscd,amssymb,verbatim,11pt]{amsart}

\theoremstyle{plain}
\newtheorem{Prop}{Proposition}[section]
\newtheorem{Thm}[Prop]{Theorem}
\newtheorem{Cor}[Prop]{Corollary}

\newtheorem{Lem}[Prop]{Lemma}

\newtheorem{Proof}[Prop]{}
\theoremstyle{definition}
\newtheorem{Def}[Prop]{Definition}

\theoremstyle{remark}
\newtheorem{Rem}[Prop]{Remark}

\numberwithin{equation}{section}
\begin{document}

\title[ Partitions of unity]%
   { Partitions of unity}

\author{
Jerzy Dydak
}
\date{ February 20, 2003
}


\thanks{ Research supported in part by a grant
 DMS-0072356 from the National Science Foundation. The author is grateful to the referee for his/her
detailed report which improved the exposition.}

\begin{abstract}

The paper contains an exposition of part of topology
using partitions of unity.
The main idea is to create variants of the Tietze
Extension Theorem and use them to derive classical theorems. 
This idea leads to a new result generalizing major
results on paracompactness (Stone Theorem and Tamano Theorem), 
a result which serves as a connection
to Ascoli Theorem. A new calculus of partitions of unity
is introduced with applications to dimension
theory and metric simplicial complexes.
The geometric interpretation of this calculus is the
barycentric subdivision of simplicial complexes.
Also, joins of partitions of unity are often used;
 they are an algebraic version of joins of simplicial
complexes.

\end{abstract}

\maketitle

\medskip
\medskip
\tableofcontents

\section{Introduction. }

The explosion of research in topology makes it imperative that one ought
to look at its foundations and decide what topics should be included in its mainstream.
One of the primary criteria is interconnectedness and potential
applications to many branches of topology and mathematics.
The author believes that the gems of basic topology are:
normality, compactness, paracompactness, and Tietze Extension Theorem.
For a unification of results on compactness see \cite{D-F}.
This paper is devoted to unification of normality, paracompactness, and Tietze Extension Theorem, a unification which leads to dimension theory and basic
geometric topology.
\par The favorite approach of general topologists to study spaces
is via open coverings (see \cite{Ho}). Geometric topologists, on the other hand, use continuous functions
to polyhedra. We plan to unify the two approaches by employing partitions of unity.
In analysis, partitions of unity form one of the basic tools. Also,
they are very useful in homotopy theory (see \cite{Di} and \cite{Do}).
In contrast,
traditional expositions of topology prove only existence of partitions of unity
subordinated to a given cover (see \cite{En$_1$} or \cite{Mu}). There is an attempt
of applying partitions of unity in \cite{L-W}. However, in \cite{L-W}
(as well as in \cite{Mu}) attention (and the definition) is restricted to locally finite
partitions of unity, in \cite{J} point finite partitions of unity are discussed. That makes applications difficult as
it is hard to construct locally finite partitions of unity 
using algebraic methods. Even arbitrary partitions of unity
form a framework too narrow to avoid all obstacles.
It turns out that equicontinuous partitions of arbitrary functions
are at the right level of  to take full advantage of calculus and algebra
of partitions of unity.
\par 
The main feature of our approach is that most of
the results follow from variations of the following
classical theorems.


\begin{Thm}[Tietze Extension Theorem for normal spaces] \label{XXX0.1} 
If $X$ is normal and $A$ is closed in $X$, then any
continuous $f:A\to [0,1]$ extends over $X$.

\end{Thm}


\begin{Thm}[(Tietze Extension Theorem for paracompact spaces] \label{XXX0.2} 
If $X$ is paracompact and $A$ is closed in $X$, then any
continuous $f:A\to E$ from $A$ to a Banach space $E$ extends over $X$.

\end{Thm} 

\ref{XXX0.1} is proved in \cite{En$_1$} (Theorem 2.1.8) and \cite{Mu} (p.219).
\ref{XXX0.2} follows from Theorem \ref{XXX4.1} in \cite{Hu}. We will subsequently
outline proofs of \ref{XXX0.1} and
\ref{XXX0.2} using our definitions of normal and paracompact spaces
(see \ref{XXXpfA} and
\ref{XXXpfB}).
\par Traditionally, topologists define a class of spaces
by using the weakest property characterizing that class.
The basic example is that of normal spaces; the definition is
that they are Hausdorff and any pair of disjoint, closed subsets
can be engulfed by disjoint, open subsets. One then
proves Urysohn Lemma and \ref{XXX0.1}-\ref{XXX0.2}.
We plan to choose one of the strongest properties
characterizing a particular class and we obtain that way
the following sequence of definitions and results which demonstrates a natural
progression of ideas (notice that some of the proofs
are postponed until subsequent sections of the paper).


\begin{Def} \label{XXX0.3} A Hausdorff space $X$ is {\bf normal}
if for any finite open covering $\{U_s\}_{s\in S}$
there is a partition of unity $\{f_s\}_{s\in S}$
on $X$ such that $f_s(X-U_s)\subseteq\{0\}$ for each $s\in S$. 

\end{Def} 

Notice that in the case of $S$ being a two-element set,
\ref{XXX0.3} is virtually identical to the Urysohn Lemma. By induction,
Urysohn Lemma implies \ref{XXX0.3}, thus showing that \ref{XXX0.3} is equivalent
to the traditional definition.
\par In this paper we shall prove several classical theorems using definitions
in our sense (\ref{XXX0.3}, \ref{XXX0.5}, and \ref{XXX0.10}). Let us start
with a natural proof of \ref{XXX0.1} using \ref{XXX0.3}.


\begin{Proof} \label{XXXpfA} Proof of \ref{XXX0.1}.
\end{Proof}
\begin{pf} The first step is to show that
any continuous function $f:A\to [-M,M]$, $A$ closed subset of a
normal (as in \ref{XXX0.3}) space $X$ extends approximately over $X$.
That is, for any $n>0$, there is a continuous function
$f_n:X\to [-M,M]$ such that $|f_n(x)-f(x)|<1/n$ for all $x\in A$.
Cover $[-M,M]$ with finitely many open (in $[-M,M]$) intervals
$\{I_s\}_{s\in S}$ of length smaller that $1/n$. 
Put $U_s=f^{-1}(I_s)\cup (X\setminus A)$ for $s\in S$
and notice that $\{U_s\}_{s\in S}$ is a finite open cover of $X$.
Pick a partition of unity $\{f_s\}_{s\in S}$ on $X$ such that
$f_s(X\setminus U_s)\subseteq \{0\}$ for all $s\in S$.
Pick $v_s\in I_s$ for each $s\in S$ and define $f_n$ via
$f_n(x)=\sum\limits_{s\in S}f_s(x)\cdot v_s$ for $x\in X$.
Notice that $|f_n(x)|\leq \sum\limits_{s\in S}f_s(x)\cdot |v_s|
\leq \sum\limits_{s\in S}f_s(x)\cdot M=M$, so the image of $f_n$ 
is in $[-M,M]$. Also, if $x\in A$, then $f_n(x)-f(x)=
\sum\limits_{s\in S}f_s(x)\cdot (v_s-f(x))$, and
$f_s(x)>0$ implies $f(x)\in I_s$ and $|v_s-f(x)|<1/n$.
Therefore, $|f(x)-f_n(x)|< \sum\limits_{s\in S}f_s(x)\cdot 1/n=1/n$
proving that $f_n$ is an approximate extension of $f$.
\par Now, given any $f:A\to [0,1]$, we construct by induction on $n$
a sequence of continuous functions $f_n:X\to [-1/2^n,1/2^n]$, $n\ge 0$,
so that $f_0=0$ and $f_{n+1}$ approximates $f-\sum\limits_{i=0}^nf_i$
within $1/2^{n+2}$. Clearly, $\sum\limits_{i=0}^\infty f_i$
is continuous and extends $f$. This extension can be modified
using a retraction $R\to [0,1]$ to get an extension of $f$
from $X$ to $[0,1]$.
\end{pf}

The main result for normal spaces ought to be one which helps proving that
certain spaces are normal. Typically, one builds spaces
from pieces, so the natural result is the one which allows
extensions of partitions of unity.


\begin{Thm} \label{XXX0.4} Suppose $X$ is normal,
$A$ is a closed subset of $X$, and $\{U_s\}_{s\in S}$
is a finite open covering on $X$.
For any finite partition of unity $\{f_s\}_{s\in S}$
on $A$ such that $f_s(A-U_s)\subseteq\{0\}$ for each $s\in S$,
there is an extension $\{g_s\}_{s\in S}$ of $\{f_s\}_{s\in S}$
over $X$
such that $g_s(X-U_s)\subseteq\{0\}$ for each $s\in S$. 

\end{Thm}

\begin{pf}  Using \ref{XXX0.1} extend each $f_s$ to $h_s:X\to [0,1]$
so that $h_s(X\setminus U_s)\subseteq \{0\}$ for $s\in S$. 
Find a neighborhood $U$ of $A$ in $X$ such that
$h=\sum\limits_{s\in S} h_s$ is positive on $U$.
Pick a continuous function $u:X\to [0,1]$
with $u(A)\subset \{1\}$ and $u(X\setminus U)\subseteq \{0\}$.
Choose a partition of unity $\{p_s\}_{s\in S}$ on $X$
so that $p_s(X\setminus U_s)\subseteq \{0\}$ for $s\in S$.
Set $q_s=u\cdot h_s+(1-u)\cdot p_s$
and notice that $q=\sum\limits_{s\in S} q_s$ is positive.
Finally, $g_s:=q_s/q$ induces the required partition of unity.
\end{pf}
 
The transition to paracompact spaces is very simple now.
Again, we do not choose the weakest condition
characterizing paracompact spaces but one of the strongest.


\begin{Def} \label{XXX0.5} A Hausdorff space $X$ is {\bf paracompact}
if for any open covering $\{U_s\}_{s\in S}$ on $X$
there is a partition of unity $\{f_s\}_{s\in S}$
on $X$ such that $f_s(X-U_s)\subseteq\{0\}$ for each $s\in S$. 

\end{Def}

Let us demonstrate how our proof of \ref{XXX0.1} can be adjusted
to give a proof of \ref{XXX0.2} using \ref{XXX0.5}.


\begin{Proof} \label{XXXpfB} Proof of \ref{XXX0.2}.
\end{Proof}
\begin{pf}  Let $E$ be a Banach space. By $B(0,M)$
we denote the open ball of radius $M$ ($M$ could be infinity) centered at $0$.
The first step is to show that
any continuous function $f:A\to B(0,M)$, $A$ being a closed subset of a
paracompact (as in \ref{XXX0.5}) space $X$,
 extends approximately over $X$.
That is, for any $n>0$, there is a continuous function
$f_n:X\to B(0,M)$ such that $|f_n(x)-f(x)|<1/n$ for all $x\in A$.
Cover $B(0,M)$ with open balls
$\{I_s\}_{s\in S}$ of diameter smaller that $1/n$. 
Put $U_s=f^{-1}(I_s)\cup (X\setminus A)$ for $s\in S$
and notice that $\{U_s\}_{s\in S}$ is an open cover of $X$.
Pick a locally finite partition of unity $\{f_s\}_{s\in S}$ on $X$ such that
$f_s(X\setminus U_s)\subseteq \{0\}$ for all $s\in S$
(their existence is shown in \ref{XXX1.12}).
Pick $v_s\in I_s$ for each $s\in S$ and define $f_n$ via
$f_n(x)=\sum\limits_{s\in S}f_s(x)\cdot v_s$ for $x\in X$.
Notice that $|f_n(x)|\leq \sum\limits_{s\in S}f_s(x)\cdot |v_s|
< \sum\limits_{s\in S}f_s(x)\cdot M=M$, so the image of $f_n$ 
is in $B(0,M)$. Also, if $x\in A$, then $f_n(x)-f(x)=
\sum\limits_{s\in S}f_s(x)\cdot (v_s-f(x))$, and
$f_s(x)>0$ implies $f(x)\in I_s$ and $|v_s-f(x)|<1/n$.
Therefore, $|f(x)-f_n(x)|< \sum\limits_{s\in S}f_s(x)\cdot 1/n=1/n$
proving that $f_n$ is an approximate extension of $f$.
\par Now, given any $f:A\to E$, let $f_0:X\to E$
approximate $f$ within $1$, i.e. $|f-f_0|<1$
and $f-f_0$ maps $A$ to $B(0,1)$. Construct by induction on $n$
a sequence of continuous functions $f_n:X\to B(0,1/2^n)$, $n\ge 1$,
so that $f_{n+1}$ approximates $f-\sum\limits_{i=0}^nf_i$
within $1/2^{n+2}$. Clearly, $\sum\limits_{i=0}^\infty f_i$
is continuous and extends $f$.
\end{pf}

Let us show an important application of \ref{XXX0.4}.


\begin{Thm} \label{XXX0.6} CW complexes are paracompact. 

\end{Thm}

\begin{pf}  Let $\{U_s\}_{s\in S}$ be an open covering on 
a CW complex $K$. For each open cell $e$ of $K$
let $K_e$ be the smallest subcomplex of $K$ containing
$e$. $K_e$ is finite for each $e$.
Our plan is to create by induction on $\dim(e)$
a partition of unity $\{f^e_s\}_{s\in S}$ on $K_e$
satisfying the following properties:
\par{a.} $f^e_s(K_e-U_s)\subseteq \{0\}$ for each $s\in S$,
\par{b.} $f^e_s\equiv 0$ for all but finitely many $s\in S$,
\par{c.} if $c\subseteq K_e$, then $f^e_s| K_c=f^c_s$
for each $s\in S$.
\par \noindent For $0$-cells $e$ it suffices to pick
one $s(e)\in S$ so that $U_{s(e)}$ contains $e$,
declare $f^e_{s(e)}\equiv 1$, and declare $f^e_s\equiv 0$
for $s\ne s(e)$.
Suppose $\{f^e_s\}_{s\in S}$ exists for all open
cells $e$ of dimension less than $n$.
Given an open $n$-cell $c$ of $K$, one can paste
all $\{f^e_s\}_{s\in S}$, $e\subseteq K_c$ and $\dim(e)<n$,
which produces a partition of unity $\{g_s\}_{s\in S}$ 
on the $(n-1)$-skeleton $L$ of $K_e$
such that $g_s\equiv 0$ for all but finitely
many $s\in S$, and $g_s(L-U_s)\subseteq \{0\}$ for each $s\in S$.
By \ref{XXX0.4} that partition of unity can be extended over $K_e$
producing $\{f^e_s\}_{s\in S}$ satisfying conditions a)-c).
\par\noindent Finally, all $\{f^e_s\}_{s\in S}$ can be pasted
together resulting in a partition of unity $\{f_s\}_{s\in S}$
on $K$ so that $f_s(K-U_s)\subseteq \{0\}$ for each $s\in S$.
\end{pf}
 
\ref{XXX0.6} is proved in \cite{L-W} (Theorem 4.2 on p.54) using Zorn Lemma
and a version of \ref{XXX0.4} for compact spaces and locally finite partitions of unity
(note that all partitions of unity in \cite{L-W} 
are assumed to be locally finite - see 
the definition on p.201 and the proof of Lemma 4.1 on p.54).
\par
Here is another illustration how a strong definition
allows for easy proofs. It also shows the advantage of using
arbitrary partitions of unity rather than only locally finite ones.


\begin{Cor} \label{XXX0.7} If $A_n$ is a closed subset of a paracompact 
(respectively, normal) space $X$
for $n\ge 1$, then $\bigcup\limits_{n=1}^\infty A_n$ is paracompact (respectively, normal). 

\end{Cor}

\begin{pf}  Let $Y=\bigcup\limits_{n=1}^\infty A_n$. Clearly, $Y$ is Hausdorff.
Suppose $\{U_s\}_{s\in S}$ is an open cover (respectively, a finite open cover) of $Y$.
Enlarge each $U_s$ to an open subset $V_s$ of $X$ so that $U_s=V_s\cap Y$. 
Notice that each $\cal V_n:=\{V_s\}_{s\in S}\cup\{X\setminus A_n\}$
is an open cover of $X$.
Pick a partition of unity $\{f_{n,s}\}_{s\in S}\cup \{f_n\}$ on $X$
for that cover.
Define $g_s:=\sum\limits_{n=1}^\infty f_{n,s}/2^n$
and $g_0:=\sum\limits_{n=1}^\infty f_{n}/2^n$.
Clearly,  $\{g_{s}\}_{s\in S}\cup \{g_0\}$ is a partition of unity
on $X$. If $g_0(x)=1$, then $f_n(x)=1$ for all $n$
and $x\in X\setminus A_n$ for all $n$ which means $x\in X\setminus Y$.
Therefore $g:=1-g_0=\sum\limits_{s\in S}g_s$ is positive on $Y$
and $h_s:=(g_{s}/g)|Y$ defines a partition of unity $\{h_{s}\}_{s\in S}$
on $Y$ such that $h_s(Y\setminus U_s)\subseteq\{0\}$ for each $s\in S$.
\end{pf}
 
In traditional approaches, \ref{XXX0.7} is much more difficult to prove
as one needs to deal first with $\sigma$-locally finite covers
(see \cite{En$_1$}, Theorem 5.1.28).
\par
Again, our main result for paracompact
spaces ought to be one which allows to extend
partitions of unity.


\begin{Thm} \label{XXX0.8} Suppose $X$ is paracompact,
$A$ is a closed subset of $X$, and $\{U_s\}_{s\in S}$
is an open covering on $X$.
For any partition of unity $\{f_s\}_{s\in S}$
on $A$ such that $f_s(A-U_s)\subseteq\{0\}$ for each $s\in S$,
there is an extension $\{g_s\}_{s\in S}$ of $\{f_s\}_{s\in S}$
over $X$
such that $f_s(X-U_s)\subseteq\{0\}$ for each $s\in S$. 
If $\{f_s\}_{s\in S}$ is locally finite, then we may
require $\{g_s\}_{s\in S}$ to be locally finite, too.
If $\{f_s\}_{s\in S}$ is point finite and $A$ is a $G_\delta$-subset
of $X$, then we may
require $\{g_s\}_{s\in S}$ to be point finite, too.

\end{Thm}

Notice that the condition of $A$ being a $G_\delta$-subset
of $X$ cannot be removed in the point finite case as shown in \cite{D$_1$}.

One can generalize the concept of the order of open covers
to partitions of unity as follows.


\begin{Def} \label{XXX0.9} A partition of unity $\{f_s\}_{s\in S}$
on $X$ is of {\bf order} at most $n$
if, for each $x\in X$, the cardinality 
of the set $\{s\in S\mid f_s(x)>0\}$ is at most $n+1$. 

\end{Def}

Now, the definition of covering dimension transfers naturally.


\begin{Def} \label{XXX0.10} A Hausdorff space $X$ is {\bf of dimension
at most $n$}
if for any open covering $\{U_s\}_{s\in S}$
there is a partition of unity $\{f_s\}_{s\in S}$
on $X$ of order at most $n$ such that $f_s(X-U_s)\subseteq\{0\}$ for each $s\in S$. 

\end{Def} 

Our main result in the theory of dimension of paracompact spaces
is the following generalization of the Tietze Extension
Theorem.


\begin{Thm} \label{XXX0.11} Let $n\ge 0$. Suppose $X$ is a paracompact space,
$\{U_s\}_{s\in S}$ is an open cover of $X$,
$A$ is a closed subset of $X$, and $\{f_s\}_{s\in S}$ is a partition of unity
on $A$ of order at most $n$ such that $f_s(A-U_s)\subseteq \{0\}$ for each $s\in S$.
There is
 a partition of unity $\{g_s\}_{s\in S}$ 
on $X$ and a closed neighborhood $B$ of $A$ in $X$ 
such that the following conditions are satisfied:
\par{a.} $g_s|A=f_s$ for each $s\in S$.
\par{b.} $g_s(X-U_s)\subseteq \{0\}$ for each $s\in S$.
\par{c.} The order of $\{g_s|B\}_{s\in S}$ is at most $n$.
\par{d.} If $\dim(X)\leq n$, then we may require $B=X$. 

\end{Thm}

Let us demonstrate the strength of our unification scheme
by discussing adjunction spaces $X\cup_fY$.
In practical applications it is important to know
that operation of taking adjunction preserves a particular
class. We will do it for normality, paracompactness,
and finite covering dimension. The feature which we would like
to emphasize is
that the proofs either change very little or form a natural
progression. Notice that even in well-known textbooks
(see \cite{Hu}, p.15) such results are only stated and their
proofs are referred to specialized papers.
\par


\begin{Def} \label{XXX0.12} Suppose $A$ is a closed
subset of a space $X$ and $f:A\to Y$ is a continuous
function. The {\bf adjunction space} $X\cup_f Y$
is the quotient space of the disjoint union
$X\oplus Y$ of $X$ and $Y$ under identification
$x\sim f(x)$ for all $x\in A$.

\end{Def}


\begin{Prop} \label{XXX0.13} Suppose $A$ is a closed
subset of a space $X$ and $f:A\to Y$ is a continuous
function. If $X$ and $Y$ are normal, then $X\cup_f Y$
is normal as well.

\end{Prop}

\begin{pf}  Suppose $B$ is a closed subset of $X\cup_f Y$
and $g:B\to I$ is a continuous function. We plan to show that $g$
extends over $X\cup_f Y$. Since $Y$ is a closed subset
of $X\cup_f Y$ and is normal, $g|B\cap Y$ extends over $Y$.
Assume then that $Y\subset B$. Passing to $X\oplus Y$,
one gets a closed subset $B'$ containing $A\oplus Y$
and a continuous function $g':B'\to I$. $g'$ can be extended over $X$,
as $X$ is normal, resulting in an extension $X\oplus Y\to I$
of $g'$. That extension induces an extension
$X\cup_f Y\to I$ of $g$.
\par Notice that one-point sets are closed in $X\cup_f Y$.
The above argument, applied to $B$ consisting of two points,
shows that for every two distinct points $x,y\in X\cup_f Y$
there is a continuous function $g:X\cup_f Y\to I$ such that $g(x)=0$
and $g(y)=1$. In particular, $X\cup_f Y$ is Hausdorff
which completes the proof.
\end{pf}


\begin{Prop} \label{XXX0.14} Suppose $A$ is a closed
subset of a space $X$ and $f:A\to Y$ is a continuous
function. If $X$ and $Y$ are paracompact, then $X\cup_f Y$
is paracompact as well. Moreover, if $\dim(X)\leq n$ and $\dim(Y)\leq n$,
then $\dim(X\cup_fY)\leq n$.

\end{Prop}

\begin{pf}  Suppose $\cal U=\{U_s\}_{s\in S}$ is an open cover of $X\cup_f Y$.
 Since $Y$ is a closed subset
of $X\cup_f Y$ and is paracompact, 
there is a partition of unity $\{g_s\}_{s\in S}$
on $Y$ such that $g_s(Y\setminus U_s)\subseteq \{0\}$
for each $s\in S$.
Let $\pi:X\oplus Y\to X\cup_f Y$ be the projection.
One gets an open cover 
$\cal V=\{V_s\}_{s\in S}$, $V_s:=\pi^{-1}(U_s)$,
and a partition of unity $\{g_s\circ \pi\}_{s\in S}$
on $A\oplus Y$ so that $(g_s\circ\pi)(A\oplus Y\setminus V_s)\subseteq \{0\}$
for each $s\in S$. By \ref{XXX0.8}, this partition can be extended over $X\oplus Y$.
 That extension induces a a partition of unity $\{h_s\}_{s\in S}$
on
$X\cup_f Y$ such that $h_s(X\cup_f Y\setminus U_s)\subseteq \{0\}$
for each $s\in S$.
\par The proof in case of both $X$ and $Y$ being of dimension
at most $n$ is exactly the same using \ref{XXX0.11}.
\end{pf}

Partitions of unity provide a simple criterion for metrizability
of a space $X$ (see \cite{D$_1$}). That criterion will be used later on
to derive the classical Bing-Nagata-Smirnov metrization theorem.


\begin{Thm} \label{XXX0.15} A Hausdorff space $X$
is metrizable if and only if there is
a partition of unity $\{f_s\}_{s\in S}$ on $X$
such that $\{f_s^{-1}(0,1]\}_{s\in S}$ is a basis
of open sets of $X$.

\end{Thm}

Notice that the classical definition of compact Hausdorff spaces
is equivalent to the following one.


\begin{Def} \label{XXX0.16} A Hausdorff space $X$ is {\bf compact}
if for any open covering $\{U_s\}_{s\in S}$
there is a finite partition of unity $\{f_s\}_{s\in S}$
on $X$ such that $f_s(X-U_s)\subseteq\{0\}$ for each $s\in S$. 

\end{Def} 

 The author does not know of any advantage in defining compact spaces that way.
However, it is useful to know that all major concepts of basic topology
can be connected using partitions of unity.

\section{Partitions of unity and equicontinuity. }


\begin{Def} \label{XXX1.1} Suppose $\{f_s\}_{s\in S}$ is a family of functions from a space $X$
to $[0,\infty)$. $\sum\limits_{s\in S}f_s=f$ means that, for each $x\in X$, $f(x)$ is the supremum of
the set $$\{\sum\limits_{s\in T}f_s(x)\mid T\text{ a finite subset of }S\}.$$
Notice that we do allow the values of $f$ to be infinity.

\end{Def} 

We are interested in families of continuous functions $\{f_s\}_{s\in S}$ from a space $X$
to $[0,\infty)$.


\begin{Def} \label{XXX1.2} A family of functions $\cal F=\{f_s:X\to [0,\infty)\}_{s\in S}$
is a {\bf partition of} 
a function $f:X\to [0,\infty]$ if $f_s$ is continuous
for each $s\in S$,
and $\sum\limits_{s\in S}f_s=f$.
In particular, $\cal F$ is called a {\bf partition of unity} 
on $X$ if $\sum\limits_{s\in S}f_s=1$.
\par\noindent
 $\cal F$ is called  a {\bf finite partition}  of $f$
 provided
$f _s\equiv 0$ for all but finitely
many $s\in S$. 
\par\noindent
$\cal F$ is called a {\bf point finite partition}  of $f$
provided $\cal F|\{x\}$ is
a  finite partition of $f|\{x\}$ for all $x\in X$.
\par\noindent
$\cal F$ is called a {\bf locally finite partition}  of $f$
provided for each $x\in X$ there is a neighborhood
$U$ of $x$ in $X$ so that $\cal F|U$ is a finite partition of $f|U$. 
\par\noindent

\end{Def}

A size of a partition $\cal F=\{f_s:X\to [0,\infty)\}_{s\in S}$ of $f$
is measured by open covers $\cal U$ of $X$ indexed by the same set $S$.


\begin{Def} \label{XXX1.3} Suppose $\cal F=\{f_s:X\to [0,\infty)\}_{s\in S}$ 
is a partition of $f:X\to [0,\infty]$ and $\cal U=\{U_s\}_{s\in S}$
is an open cover of $X$.
$\cal F$ is {\bf $\cal U$-small} if 
 $f_s(X-U_s)\subseteq\{0\}$ for each $s\in S$. In other words,
the carrier of $f_s$ is contained in $U_s$ for each $s\in S$.

\end{Def}

The goal of this section is to produce partitions of
continuous functions. We are going to proceed in small, simple steps.
The first order of business is to characterize continuity of $\sum\limits_{s\in S}f_s$
in terms similar to those for power series. In our case
a \lq tail\rq\ of $\sum\limits_{s\in S}f_s$ is
$\sum\limits_{s\in S\setminus T}f_s$, $T$ finite subset of $S$.


\begin{Prop} \label{XXX1.4} Suppose $\{f_s\}_{s\in S}$
is a family of continuous functions from a space $X$
to $[0,\infty)$ so that $\sum\limits_{s\in S}f_s=f$ is finite (i.e., 
 $f(X)\subseteq [0,\infty)$). 
$f$ is continuous if and only if for each point $x\in X$
and each $\epsilon>0$ there is a neighborhood $U$
of $x$ in $X$ and a finite subset $T$ of $S$
such that the values of $\sum\limits_{s\in S\setminus T}f_s$
on $U$ are less that $\epsilon$.

\end{Prop}

\begin{pf}  For any finite subset $T$
of $S$ let $f_T$ be defined as $\sum\limits_{s\in T}f_s$.
 If $f$ is continuous, $\epsilon>0$, and $x\in X$,
then we pick a finite $T\subseteq S$ such that $f(x)-f_T(x)<\epsilon/3$.
Since $f-f_T$ is continuous, there is a neighborhood $U$ of $x$
such that $f(y)-f_T(y)<\epsilon$ for all $y\in U$.
Since $f-f_T=\sum\limits_{s\in S\setminus T}f_s$, we are done with the first
implication.
\par Suppose $U$ is a neighborhood
of $x$ in $X$ and $T$ is a finite subset of $S$
such that the values of $f-f_T=\sum\limits_{s\in S\setminus T}f_s$
on $U$ are less that $\epsilon/3$.
Find a neighborhood $V$ of $x$ in $U$ such that
$|f_T(y)-f_T(x)|<\epsilon/3$ for each $y\in V$.
Now, $|f(y)-f(x)|\leq |f(y)-f_T(y)|+|f_T(y)-f_T(x)|+|f_T(x)-f(x)|<\epsilon$
for all $y\in V$ which proves continuity of $f$ at $x$.
\end{pf}

\begin{Rem}  In \cite{Y$_1$} K.Yamazaki
calls a collection $\{ f_s\}_{s\in S}$ of continuous non-negative real-valued function on a topological space $X$ {\bf sum-complete} if $\sum\limits_{s\in S}f_s$ is a continuous function from $X$ into $[0, \infty)$, and proved that the property that every sum-complete collection of functions on 
a subset $A$ can be extended to a sum-complete collection of functions on $X$ is equivalent to
 $A$ being $P$-embedded in  $X$.
 
\end{Rem}


\begin{Cor} \label{XXX1.5} Suppose $\{f_s\}_{s\in S}$ and $\{g_s\}_{s\in S}$
are two families of continuous functions from a space $X$
to $[0,\infty)$ so that $\sum\limits_{s\in S}f_s=f$ is continuous
and $f(X)\subseteq [0,\infty)$. If $g_s(x)\leq f_s(x)$
for each $x\in X$ and each $s\in S$, then $g=\sum\limits_{s\in S}g_s:X\to [0,\infty)$ is
continuous. 

\end{Cor}

\begin{pf}  The tails of $\{g_s\}_{s\in S}$ are estimated from above
by the tails of $\{f_s\}_{s\in S}$. 
\end{pf}

Notice that if the tails of $\{f_s\}_{s\in S}$ are small,
then the family $\{\max(0,f_s-\epsilon)\}_{s\in S}$ is locally finite
for any $\epsilon>0$. This leads to a new concept. It implies
equicontinuity of $\{f_s\}_{s\in S}$, hence its name.


\begin{Def} \label{XXX1.6} Suppose $\{f_s\}_{s\in S}$ is a family of continuous functions 
from a space $X$
to $[0,\infty)$. $\{f_s\}_{s\in S}$ is called
{\bf strongly equicontinuous} if one of the following 
equivalent conditions holds:
\par{a.} For each $\epsilon>0$ and each $x\in X$
there is a neighborhood $U$ of $x$ in $X$ and a finite subset
$T$ of $S$ such that
 $f_s(y)<\epsilon$ for all $y\in U$ and all $s\in S\setminus T$.
\par{b.} For each positive $\epsilon$
the family
$\{\max(0,f_s-\epsilon)\}_{s\in S}$ is locally finite.

\end{Def}

Recall the concept of
equicontinuity (see 3.4.17 in \cite{En$_1$} or \cite{Mu}, p.276).


\begin{Def} \label{XXX1.7}
A family of functions $\{f_s\}_{s\in S}$
from a space $X$ to a metric space $(Y,d)$
is {\bf equicontinuous} if for each $\epsilon>0$ and each point $a\in X$ there
a neighborhood $U$ of $a$ in $X$ such that $d(f_s(x),f_s(y))<\epsilon$
for all $s\in S$ and all $x,y\in U$. 

\end{Def}

Next we show that strong equicontinuity implies equicontinuity.
Surprisingly, if the sum of functions is finite,
they are equivalent.


\begin{Prop} \label{XXX1.8} Suppose $\{f_s\}_{s\in S}$ is a family of continuous functions 
from a space $X$
to $[0,\infty)$. Consider the following conditions:
\par{a.} $\{f_s\}_{s\in S}$ is strongly equicontinuous.
\par{b.} $\{\max(0,f_s-g)\}_{s\in S}$ is locally finite
for any positive, continuous $g:X\to R$.
\par{c.}
$\{f_s\}_{s\in S}$ is equicontinuous.
\par Conditions a) and b) are equivalent. Condition a) implies Condition c).
If $\sum\limits_{s\in S}f_s=f$ is finite (i.e., 
 $f(X)\subseteq [0,\infty)$), then all three conditions are equivalent.

\end{Prop}

\begin{pf}  a)$\implies$ b). Given $x\in X$ and a positive continuous
function $g:X\to R$ put $\epsilon=g(x)/2$ and find
a neighborhood $V$ of $x$ in $X$ such that $g(y)>\epsilon$ for all $y\in V$.
By a) there is a neighborhood $U$ of $x$ in $V$ and a finite subset $T$ of $S$
such that $f_s(y)\leq \epsilon$ for all $y\in U$ and all $s\in S\setminus T$.
Notice that $\max(0,f_s(y)-g(y))=0$ for all 
$y\in U$ and all $s\in S\setminus T$ which proves b).
\par b)$\implies$ a). Suppose $\epsilon>0$. Put $g\equiv\epsilon$.
\par a)$\implies$ c). Given $\epsilon>0$ and $x\in X$
find a finite subset $T$ of $S$ and a neighborhood
$V$ of $x$ such that $f_s(y)\leq \epsilon/3$ for all $y\in V$ and all $s\in S\setminus T$.
In particular, $|f_s(z)-f_s(y)|<\epsilon$
for all $s\in S\setminus T$ and all $z,y\in V$. 
Obviously, $\{f_s\}_{s\in T}$ is equicontinuous,
so there is a neighborhood $U$ of $x$ in $V$ such that
$|f_s(z)-f_s(y)|<\epsilon$
for all $s\in T$ and all $z,y\in U$. 
\par Assume $\sum\limits_{s\in S}f_s=f$ is finite and
$\{f_s\}_{s\in T}$ is equicontinuous.
Given $x\in X$ and $\epsilon>0$ pick
a neighborhood $U$ of $x$ in $X$ such that
$|f_s(z)-f_s(y)|<\epsilon$
for all $s\in S$ and all $z,y\in U$. 
Find a finite subset $T$ of $S$ such that
$f(x)-\sum\limits_{s\in T}f_s(x)<\epsilon/2$.
That implies $f_s(x)<\epsilon/2$ for all $s\in S\setminus T$.
Now, if $s\in S\setminus T$ and $y\in U$, then
$f_s(y)<f_s(x)+\epsilon/2<\epsilon/2+\epsilon/2=\epsilon$.
\end{pf}

From now on we will be interested in producing equicontinuous families
of functions. Therefore the following simple fact is useful
as it can be applied to $\{f_s+g_t\}_{s,t\in T}$,
$\{\max(0,f_s-g_t)\}_{s,t\in T}$, and so on.


\begin{Prop} \label{XXX1.9} Suppose $\{f_s:X\to Y\}_{s\in S}$ and
$\{g_t:X\to Y'\}_{t\in T}$ are two families of
functions from a space $X$ to metric spaces $(Y,d)$ and $(Y',d')$. 
\par{a.} If $h:Y\to Y'$ is uniformly continuous
and $\{f_s:X\to Y\}_{s\in S}$ is equicontinuous,
then $\{h\circ f_s:X\to Y'\}_{s\in S}$ is equicontinuous.
\par{b.} The two families are equicontinuous
if and only if $\{h_{s,t}:X\to Y\times Y'\}_{s\in S,t\in T}$
defined by $h_{s,t}(x)=(f_s(x),g_t(x))$ is equicontinuous.

\end{Prop}

\begin{pf}  a). Suppose $a\in X$ and $\epsilon>0$.
Choose $\delta>0$ such that $d(y_1,y_2)<\delta$ implies
$d'(h(y_1),h(y_2))<\epsilon$. Let $U$ be a neighborhood
of $a$ in $X$ such that $d(f_s(x),f_s(y))<\delta$
for all $x,y\in U$. Now,
$d'(h\circ f_s(x),h\circ f_s(y))<\epsilon$ for all
$x,y\in U$.
\par b). $Y\times Y'$ is considered with the
metric $\rho$ being the sum of $d$ and $d'$.
Notice that the projections $Y\times Y'\to Y$
and $Y\times Y'\to Y'$ are uniformly continuous, so a) implies part of b).
\par
Suppose $\epsilon>0$ and $a\in X$. Find neighborhoods $U$ and $U'$
of $a$ in $X$ such that $d(f_s(x),f_s(y))<\epsilon/2$ for each $x,y\in U$
and each $s\in S$,
and $d(g_t(x),g_t(y))<\epsilon/2$ for each $x,y\in U'$ and each $t\in T$.
Notice that $\rho(h_{s,t}(x),h_{s,t}(y))<\epsilon$
for each $x,y\in U\cap U'$.
\end{pf}

Notice that the classical concepts of the supremum and the infimum of a subset of reals
can be naturally extended to the concepts of the supremum $\sup \{f_s\}_{s\in S}$
and the infimum $\inf \{f_s\}_{s\in S}$ of
a family of real-valued functions on any set $X$.

The following result is crucial in production of equicontinuous families
of functions.


\begin{Prop} \label{XXX1.10} Suppose $\{f_s\}_{s\in S}$ is an equicontinuous family of
functions from a space $X$ to reals $R$. If $\sup \{f_s\}_{s\in S}<\infty$
(respectively, $\inf \{f_s\}_{s\in S}>-\infty$),
then the family $\{f_T\}_{T\subseteq S}$ is equicontinuous,
where $f_T:=\sup \{f_s\}_{s\in T}$ (respectively, $f_T:=\inf \{f_s\}_{s\in T}$).

\end{Prop}

\begin{pf}  Suppose $a\in X$ and $\epsilon>0$. We need to find
a neighborhood $U$ of $a$ in $X$ such that $|f_T(x)-f_T(y)|<\epsilon$
for all $T\subseteq S$ and all $x,y\in U$.
Let $U$ be a neighborhood of $a$ in $X$ such that $|f_s(x)-f_s(y)|<\epsilon/2$
for all $s\in S$ and all $x,y\in U$.
 It suffices to show $f_T(x)< f_T(y)+\epsilon$
for all $T\subseteq S$ and all $x,y\in U$ (use symmetry).
Since $f_s(x)<f_s(y)+\epsilon/2\leq f_T(y)+\epsilon/2$
for all $s\in T$, taking the supremum of the left side
results in $f_T(x)\leq f_T(y)+\epsilon/2<f_T(y)+\epsilon$.
\end{pf}

The following concept will be useful.


\begin{Def} \label{XXX1.11} A partition of unity $\{g_s\}_{s\in S}$ on $X$ 
is an {\bf approximation} of
a partition $\{f_s\}_{s\in S}$ of $f$ if $g_s(x)>0$ implies $f_s(x)>0$ for every
$s\in S$.

\end{Def}


\begin{Cor} \label{XXX1.12} Every equicontinuous partition $\{f_s\}_{s\in S}$
of a positive and finite function $f:X\to (0,\infty)$
has a locally finite approximation $\{g_s\}_{s\in S}$
such that the closure of the carrier of $g_s$ is contained
in the carrier of $f_s$ for each $s\in S$. 

\end{Cor}

\begin{pf}  By replacing $f_s$ with $\min(1,f_s)$ we may assume that $f_s:X\to [0,1]$
for each $s\in S$.
Let $g:=\sup\{f_s\mid s\in S\}$ and $h_s:=\max(0,f_s-g/2)$.
$g$ is continuous by \ref{XXX1.10} and positive-valued.
Also, for each $a\in X$, there is $s\in S$ with $g(s)=f_s(a)>0$
which implies that $h_s(a)=g(s)/2>0$. By \ref{XXX1.8} functions
$h_s$ induce a $\cal U$-small, locally finite
partition of a continuous, positive-valued function $h:X\to (0,\infty)$
such that the closure of the carrier of $h_s$ is contained
in the carrier of $f_s$ for each $s\in S$.
Put $g_s:=h_s/h$.
\end{pf}

Given an open cover $\cal U=\{U_s\}_{s\in S}$ of a space $X$
it is natural to seek sufficient conditions for
existence of a $\cal U$-small partition of unity on $X$.
In case of countable covers one has a simple
necessary and sufficient condition.


\begin{Prop} \label{XXX1.13} Suppose $\cal U=\{U_n\}_{n\ge 1}$ 
is a countable open cover of a space $X$.
A $\cal U$-small partition of unity on $X$ exists if and only if
there is a positive-valued $f:X\to (0,\infty]$ which has a $\cal U$-small
partition.

\end{Prop}

\begin{pf}  Suppose $\cal F=\{f_n\}_{n\ge 1}$ 
is a $\cal U$-small partition of $f:X\to (0,\infty]$.
Put $g_n=\min(f_n,2^{-n})$ for $n\ge 1$.
It is well-known that it is a partition of
a continuous $g$. Alternatively,
 notice that the tails of $\{g_n\}_{n\ge 1}$ are small
and use \ref{XXX1.4}.
Therefore $\{g_n/g\}_{n\ge 1}$ 
is a $\cal U$-small partition of unity on $X$.
\end{pf}

\ref{XXX1.13} immediately implies that all separable metric spaces $X$
are paracompact. Indeed, one can reduce the question of existence
of partitions of unity to countable open covers
$\cal U=\{U_n\}_{n\ge 1}$ of $X$ for which $f_n(x):=dist(x,X-U_n)$, $x\in X$,
defines a $\cal U$-small partition of a positive-valued
$f:X\to (0,\infty]$.
\par For arbitrary metric spaces one has to work with
the family of $f_s(x):=dist(x,X-U_s)$. That family does not have
to be a partition of a continuous $f:X\to [0,\infty)$
but it has an important property (see \ref{XXX2.1}) of such partitions.
\par The following is our weakest condition
characterizing existence of a $\cal U$-small
partition of unity. We will see later that it 
implies all major theorems on paracompactness.


\begin{Thm} \label{XXX1.14} Suppose $\cal U=\{U_s\}_{s\in S}$ 
is an open cover of a normal space $X$.
A $\cal U$-small partition of unity on $X$ exists if and only if
there is an equicontinuous
family $\{f_t\}_{t\in T}$ 
satisfying the following two conditions:
\par{1.} For each $x\in X$ there is $t\in T$
so that $f_t(x)>0$.
\par{2.} For each $t\in T$ there is a finite subset $F$ of $S$
with the property that $f_t(x)=0$ for 
all $x\in X\setminus\bigcup\limits_{s\in F}U_s$.

\end{Thm}

The proof of \ref{XXX1.14} is preceeded by a lemma. The purpose of 
it is to create a $\cal U$-small equicontinuous partition of a bounded function
so that we can use \ref{XXX1.12}.


\begin{Lem} \label{XXX1.15} 
Suppose $\cal U=\{U_s\}_{s\in S}$ 
is an open cover of a space $X$ and $\{f_s:X\to [0,1]\}_{s\in S}$
is a $\cal U$-small, equicontinuous partition of a positive-valued $f:X\to (0,\infty]$.
 If $S$ is well-ordered, then functions
$g_s:=\max(0,f_s-\sup\{f_t\mid t<s\})$
induce a $\cal U$-small, equicontinuous partition $\{g_s\}_{s\in S}$
of a positive-valued $g:X\to (0,1]$.

\end{Lem}

\begin{pf}  Given $a\in X$ let $t$ be the smallest element of
$\{s\in S\mid f_s(a)>0\}$. Since $g_t(a)=f_t(a)$, it follows
that $\sum\limits_{s\in S}g_s>0$.  
\par Suppose $T$ is a finite subset of $S$ such that
$g_s(a)>0$ for each $s\in T$. Enumerate all elements of $T$
in the increasing order $s(1)<s(2)<\ldots<s(k)$.
Now, $g_{s(i)}(a)\leq f_{s(i)}(a)-f_{s(i-1)}(a)$ for $i=2,\ldots,k$,
so $$\sum\limits_{s\in T}g_s(a)\leq f_{s(1)}(a)+(f_{s(2)}(a)-f_{s(1)}(a))+\ldots+
(f_{s(k)}(a)-f_{s(k-1)}(a))=f_{s(k)}(a)\leq 1$$
which proves $g(a)=\sum\limits_{s\in S}g_s(a)\leq 1$.
\par The equicontinuity of $\{g_s\}_{s\in S}$ follows from
\ref{XXX1.9}. Indeed, $\max(u,v)=(u+v+|u+v|)/2$ for any $u,v\in R$.
\end{pf}


\begin{Proof} \label{XXXpfC} Proof of \ref{XXX1.14}.
\end{Proof}
\begin{pf}
Suppose $\cal U=\{U_s\}_{s\in S}$ 
is an open cover of a normal space $X$. Pick an equicontinuous
partition $\{f_t\}_{t\in T}$ of a positive-valued function $f$ so that
for any $t\in T$ there is a finite subset $F(t)$ of $S$
with the property that $f_t(x)=0$ for 
all $x\in X\setminus\bigcup\limits_{s\in F(t)}U_s$.
Replacing $f_t$ by $\min(1,f_t)$ and using \ref{XXX1.9}, \ref{XXX1.15},
we may assume that $f$ is bounded by $1$.
By \ref{XXX1.12} we may assume $\{f_t\}_{t\in T}$ is a locally finite
partition of unity.
For each finite $F\subseteq S$ define $U_F= \bigcup\limits_{s\in F}U_s$
and let $f_T$ be the sum of all $f_t$
such that $F=F(t)$.
Clearly $\{f_F\}_{F\subseteq S}$ is a
 $\cal U'$-small partition of unity, where $\cal U'=\{U_F\}_{F\subseteq S}$.
By \ref{XXX1.12} there is a locally finite partition of
unity $\{g_F\}_{F\subseteq S}$ with the property
that the closure $A_F$ of the carrier of $g_F$ is contained in $U_F$
for each $F$.
Given $F$ consider the open cover $\{A_F\cap U_s\}_{s\in F}$ of $A_F$
and pick a partition of unity $\{h_{F,s}\}_{s\in F}$
on $A_F$ so that $h_{F,s}(A_F\setminus U_s)\subseteq \{0\}$ for each
$s\in F$ (see \ref{XXX0.4}).
We can extend each $h_{F,s}$ over $X$ so that $h_{F,s}(X-U_s)\subseteq\{0\}$.
Notice that $h_{F,s}\cdot g_F$ with $s\in F$ and $F$ ranging over all
finite subsets of $S$ forms a partition of unity on $X$.
Therefore $p_s:=\sum\limits_{F\subseteq S}h_{F,s}\cdot g_F$
induces a partition of unity on $X$. Clearly, it is $\cal U$-small.
\end{pf}

\section{Applications to general topology. }


\begin{Lem} \label{XXX2.1} Suppose $X$ is a metric space
 and $\cal U=\{U_s\}_{s\in S}$ is a family of
open subsets in $X$. The family $\cal F=\{f_s\}_{s\in S}$
of functions defined by $f_s(x):=dist(x,X-U_s)$
is equicontinuous.

\end{Lem}

\begin{pf}  For each $z\in X$ define
$g_z:X\to R$ by $g_z(x)=d(x,z)$. The Triangle Inequality
implies $|g_z(x)-g_z(y)|\leq d(x,y)$ for all
$x,y,z\in X$. In particular, $\{g_z\}_{z\in X}$
is equicontinuous. By \ref{XXX1.10},
the family $\{g_T\}_{T\subseteq X}$ is equicontinuous,
where $g_T:=\inf\{g_z\mid z\in T\}$. 
Taking $T=X-U_s$ gives $g_T=f_s$ and completes the proof.
\end{pf}

\ref{XXX1.14} and \ref{XXX2.1} imply the famous theorem of A.H.Stone (see \cite{En$_1$}, 4.4.1 and 5.1.3).


\begin{Cor}[A.H.Stone] \label{XXX2.2} 
Every metrizable space $X$
is paracompact. 

\end{Cor}

The following result describes a useful family of
 equicontinuous functions. It is well-known but the proof is so short
that we include it. In the Appendix we will show
that all equicontinuous families with values in compact spaces
are detected that way and we will apply it to prove
a basic version of Ascoli Theorem.


\begin{Lem} \label{XXX2.3} If $Z$ is a metric space
and $Y$ is a compact space, then any
continuous function $f:X\times Y\to Z$
induces an equicontinuous family
$\{f_y:X\to Z\}_{y\in Y}$ given by $f_y(x)=f(x,y)$
for $y\in Y$ and $x\in X$.

\end{Lem}

\begin{pf}  Given $a\in X$ and $\epsilon>0$ one can find,
for each $y\in Y$ a neighborhood $U_y\times V_y$
of $(a,y)$ in $X\times Y$ such that
$d(f(u,v),f(a,y))<\epsilon/2$ for all 
$(u,v)\in U_y\times V_y$.
Pick a finite cover $V_{y(1)}\cup\ldots\cup V_{y(k)}$ of $Y$
and put $U=\bigcap\limits_{i=1}^kU_{y(i)}$.
\end{pf}


\begin{Thm}[H.Tamano \cite{En$_1$}, 5.1.38] \label{XXX2.4} 
Suppose $X$ is a completely regular space. If $X\times rX$
is normal for some compactification $rX$ of $X$,
then $X$ is paracompact. 

\end{Thm}

\begin{pf}  Suppose $\cal U=\{U_s\}_{s\in S}$ is an open cover of $X$.
Obviously, if $\cal U$ has a finite subcover,
then there is a $\cal U$-small partition of unity on $X$ (see \ref{XXX0.4}).
Therefore we assume that $\cal U$ has no finite subcover.
Enlarge each $U_s$ to an open subset $V_s$ or $rX$
satisfying $X\cap V_s=U_s$.
Let $C=rX\setminus\bigcup\limits_{s\in S}V_s$
and let $A=\{(x,x)\mid x\in X\}\subseteq X\times rX$.
Notice that $C$ is non-empty, $A$ does not intersect $X\times C$,
and $A$ is closed (it is the intersection of the diagonal
in $rX\times rX$ with $X\times rX$).
Choose a continuous $f:X\times rX\to [0,1]$
so that $f(A)=\{0\}$ and $f(X\times C)=\{1\}$.
The functions $\{f_z:X\to [0,1]\}_{z\in rX}$ defined
by $f_z(x)=f(x,z)$ form an equicontinuous
family by \ref{XXX2.3}. Therefore (see \ref{XXX1.10})
$f_T:=\inf\{f_z\mid z\in X\setminus \bigcup\limits_{s\in T}U_s\}$
form an equicontinuous family of functions,
where $T$ ranges over finite subsets of $S$.
Since $f_x(x)=0$ for each $x\in X$, we get $f_T(X\setminus \bigcup\limits_{s\in T}U_s)\subset \{0\}$
for each $T\subseteq S$. It remains to show that $\sum\limits_{T\subseteq S}f_T$
is positive in view of \ref{XXX1.14}.
Given $x\in X$ find a neighborhood $U$ of $C$ in $rX$ 
with the property that $f(x,z)>1/2$ for each $z\in U$.
$rX\setminus U$ is compact and contained in $\bigcup\limits_{s\in S}V_s$,
so there is a finite subset $T$ of $S$ with the property
$rX\setminus U\subseteq \bigcup\limits_{s\in T}V_s$.
Therefore $X-\bigcup\limits_{s\in T}U_s\subseteq U$
and $f_T(x)\ge 1/2$.
\end{pf}

Notice (see \cite{En$_1$}, 5.2.A) that a Hausdorff space $X$
is normal countably paracompact if and only if any countable open cover
$\cal U=\{U_n\}_{n\in N}$ of $X$ admits a partition of unity
$\{f_n\}_{n\in N}$ on $X$ such that $f_n(X\setminus U_n)\subseteq \{0\}$
for each $n\ge 1$.

Another corollary to \ref{XXX1.14} and \ref{XXX2.3} is 
the following sufficient condition
for countable paracompactness (see \cite{En$_1$}, 5.2.8 and 5.2.H).


\begin{Thm} \label{XXX2.5} 
Suppose $X$ is a space. If $X\times Y$
is normal for some infinite compact Hausdorff space $Y$,
then $X$ is countably paracompact. 

\end{Thm}

\begin{pf}  First consider $Y=rN$ to be a compactification of the natural numbers.
Suppose $\cal U=\{U_n\}_{n\in N}$ is an open cover of $X$.
Put $V_n=\bigcup\limits_{i=1}^nU_i$
and consider $A=\bigcup\limits_{n=1}^\infty(X-V_n)\times\{n\}$,
$B=X\times (rN\setminus N)$. $B$ is clearly closed,
and $A$ is closed as its complement is $\bigcup\limits_{n=1}^\infty V_n\times W_n$,
where $W_n=rN\setminus \{1,\ldots,n-1\}$.
Choose a continuous $f:X\times rN\to [0,1]$
so that $f(A)=\{0\}$ and $f(X\times (rN\setminus N))=\{1\}$.
The functions $\{f_n:X\to [0,1]\}_{n\in N}$ defined
by $f_n(x)=f(x,n)$ form an equicontinuous
family (see \ref{XXX2.3}). 
Obviously, $f_n(X-V_n)\subseteq\{0\}$ for each $n\in N$.
 To use \ref{XXX1.14}, it remains to show that $\sum\limits_{n\in N}f_n$
is positive.
If $f_n(x)=0$ for all $n\in N$, then $f(x,z)=0$
for all $z\in rN\setminus N$
contradicting $f(X\times (rN\setminus N))=\{1\}$.
\par To complete the proof notice that any infinite compact Hausdorff
space contains a compactification of natural numbers.
\end{pf}

Finally, we will see how to get classical
results of general topology via the \lq discretization process\rq\
of replacing partitions of unity
by closed covers.
Here is a well-known discrete interpretation of normal spaces.


\begin{Prop} \label{XXX2.6} A Hausdorff space $X$ is normal if and only if for any
finite open cover $\{U_s\}_{s\in S}$ of $X$ there is
a closed cover  $\{F_s\}_{s\in S}$ of $X$ such that
$F_s\subseteq U_s$ for each $s\in S$. 

\end{Prop}

\begin{pf}  One direction follows from \ref{XXX0.3}.
For the other implication, use \ref{XXX0.4} and \ref{XXX1.12}.  
\end{pf}

Here is the corresponding result for paracompactness.
The challenge is to demonstrate the discretization
of the proof of \ref{XXX1.14}.


\begin{Thm}[Michael \cite{M} or \cite{En$_1$}, 5.1G] \label{XXX2.7}  
 A Hausdorff space $X$ is paracompact if and only if for any
 open cover $\cal U=\{U_s\}_{s\in S}$ of $X$ there is
a closed cover  $\{F_s\}_{s\in S}$ of $X$ such that
$F_s\subseteq U_s$ for each $s\in S$,
and $\bigcup\limits_{s\in T} F_s$ is closed for every subset $T$ of $S$. 

\end{Thm}

\begin{pf}  Obviously, one direction follows quickly from \ref{XXX1.12}.
It is the other implication which is of interest.
Our proof starts as that in \cite{En$_1$}, 5.1.33 (how else?) but
 is simpler and is motivated by partitions of unity.
As in \ref{XXX1.15}, we assume that $S$ is well-ordered 
and our basic idea is to follow the recipe of 
replacing $f_s$ by $f_s-\sup\{f_t\mid t<s\}$ adapted to the discrete case.
\par We will create, for each $n\ge 1$, closed covers $\cal F_n=\{F_{s,n}\}_{s\in S}$ of $X$
which are closure-preserving (that means
$\bigcup\limits_{s\in T} F_{s,n}$ is closed for every subset $T$ of $S$)
and $F_{s,n}\subset U_s$ for each $s\in S$.
Covers $\cal F_n$ are required to have the property that 
$$F_{s,n+1}\subseteq U_s\setminus \bigcup\limits_{t<s}\bigcup\limits_{k\leq n}F_{t,k}$$
for each $s\in S$ and $n\ge 1$.
$\cal F_1$ can be chosen by our hypotheses.
Assume $\cal F_k$ exists for $k\leq n$. Notice that 
$V_s:=U_s\setminus \bigcup\limits_{t<s}\bigcup\limits_{k\leq n}F_{t,k}$
cover all of $X$ if $s$ runs through $S$. Indeed, given $x\in X$ one can find
the smallest $s\in S$ with $x\in U_s$ in which case $x\in V_s$.
Therefore, we pick a closed, closure-preserving cover $\cal F_{n+1}=\{F_{s,n+1}\}_{s\in S}$
such that $F_{s,n+1}\subseteq U_s\setminus \bigcup\limits_{t<s}\bigcup\limits_{k\leq n}F_{t,k}$
for each $s\in S$.
\par Each $x\in X$
has a natural system of neighborhoods $W_{x,k}$, where $W_{x,k}$ is
defined as the complement of all $F_{s,p}$ not containing $x$
so that $p\leq k$. $\{W_{x,k}\}_{k\ge 1}$ is our initial approximation of neighborhoods
of $x$ needed to establish equicontinuity.
\par Our first observation is that $W_{x,k}\cap F_{s,n}\ne\emptyset$
and $k>n$ implies $x\in F_{s,n}$.
The second observation is that that $W_{x,k}\cap F_{s,n}\cap F_{s,n+1}\ne\emptyset$
and $k>n+1$ implies that $s$ is the smallest
element of $\{t\in S\mid x\in F_{t,n}\}$. Indeed,
$x\in F_{s,n+1}$ means that $x$ cannot belong to $F_{t,n}$ for any $t<s$.
Notice that the second observation implies that
$\{F_{s,n}\cap F_{s,n+1}\}_{s\in S}$ is a discrete family: 
if $m>n+1$, then $W_{x,m}$ intersects at most one of those elements.
Finally, our third observation
is that $F_{s,n}\cap F_{s,n+1}\cap F_{s,n+2}
\subseteq X-\bigcup\limits_{t\ne s} F_{t,n+1}\subseteq F_{s,n+1}$.
It follows from the fact that $x\in F_{t,n+1}$, $t<s$,
implies $x\notin F_{s,n+2}$, and 
$x\in F_{t,n+1}$, $t>s$, implies $x\notin F_{s,n}$.
Now, $E_{s,n}:=F_{s,n}\cap F_{s,n+1}\cap F_{s,n+2}\cap F_{s,n+3}
\subseteq V_{s,n}:=(X-\bigcup\limits_{t\ne s} F_{t,n+1})\cap
(X-\bigcup\limits_{t\ne s} F_{t,n+2})\subseteq F_{s,n+1}\cap F_{s,n+2}$
and $\{V_{s,n}\}_{s\in S}$ is a discrete family of open sets.
\par
Given $x\in X$ find the smallest $s$
with $x\in F_{s,n}$ for some $n$. Now, $x\notin F_{t,n+k}$ if $t<s$
and $x\notin F_{t,n+k}$ if $t>s$ and $k\ge 1$ (as such $F_{t,n+k}$ is disjoint with
$F_{s,n}$), which implies $x\in E_{s,n}$. That means sets $E_{s,n}$
cover $X$. 
\par
For each $(s,n)\in S\times N$ pick a continuous function $f_{s,n}:X\to [0,1/n]$
so that $f_{s,n}(X\setminus V_{s,n})\subset \{0\}$
and $f_{s,n}(E_{s,n})\subseteq\{1/n\}$.
\par Notice that $\{f_{s,n}\}_{(s,n)\in S\times N}$ is equicontinuous.
Indeed, given $x\in X$ and $\epsilon>0$ we can find $n\in N$
with $\epsilon>1/n$. Since each family $\cal V_k:=\{V_{s,k}\}_{s\in S}$
is discrete, we can find a neighborhood $W:=W_{x,n+1}$ of $x$ in $X$
intersecting at most one element of $\cal V_k$ for $k\leq n$.
Now, the set of non-zero $\{f_{s,k}|W\}_{s\in S, k\leq n}$ is finite,
hence equicontinuous, so there is a neighborhood $U$ of $x$ in $W$
so that $y,z\in U$ implies $|f_{s,k}(y)-f_{s,k}(z)|<\epsilon$
for all $s\in S$ and $k\leq n$.
If $k>n$, then $|f_{s,k}(y)-f_{s,k}(z)|\leq 1/n<\epsilon$.
Use \ref{XXX1.14} to conclude the proof.
\end{pf}


\begin{Cor}[Michael \cite{En$_1$},5.1.33] \label{XXX2.8} 
If $f:X\to Y$ is a closed continuous function
and $X$ is paracompact, then $Y$ is paracompact. 

\end{Cor}

\begin{pf}  It follows from \ref{XXX2.7}. 
\end{pf}


\begin{Rem} \label{XXX11.1} Notice that one can easily adapt \ref{XXX2.7} to countable covers
and conclude that images under closed continuous functions
of countably paracompact, normal spaces are countably paracompact
(see \cite{En$_1$}, 5.2.G(e)). 

\end{Rem} 


\begin{Proof} \label{XXXpfD} Proof of \ref{XXX0.15}.
\end{Proof}
\begin{pf} Suppose $X$ is metrizable and $d$ is a metric on $X$.
By \ref{XXX2.2}, given $n\ge 1$ pick a partition of unity
$\{f_{x,n}\}_{x\in X}$ such that $f_{x,n}(X\setminus B(x,1/n))\subseteq \{0\}$,
where $B(x,1/n):=\{y\in X| d(x,y)<1/n\}$
is the open $(1/n)$-ball centered at $x$.
Define $g_{x,n}:=f_{x,n}/2^n$ for $(x,n)\in X\times N$
and notice that $\{g_{x,n}\}_{(x,n)\in X\times N}$
is a partition of unity on $X$ such that
$\{g_{x,n}^{-1}(0,1]\}_{(x,n)\in X\times N}$
is a basis of open neighborhoods of $X$.
\par Suppose $X$ is a Hausdorff space and
$\{f_{s}\}_{s\in S}$
is a partition of unity on $X$ such that
$\{f_{s}^{-1}(0,1]\}_{s\in S}$
is a basis of open neighborhoods of $X$.
Define $d(x,y):=\sum\limits_{s\in S}|f_s(x)-f_s(y)|$.
It is clearly a metric on $X$, so it remains to show that
it induces the same topology on $X$.
Suppose $U$ is an open set in $X$ and $x\in U$.
There is $t\in S$ such that $x\in f_t^{-1}(0,1]\subset U$.
Consider $V:=\{y\in X| d(x,y)<f_t(x)\}$.
Notice that $V$ is open and contains $x$. To show $V\subset U$
assume $y\in V\setminus U$. Now $f_t(y)=0$,
so $d(x,y)\ge |f_t(x)-f_t(y)|=f_t(x)$, a contradiction.
\end{pf}

We are ready to derive part of the classical Nagata-Smirnov metrizability
criterion \cite{En$_1$}. The second part will be derived later on (see \ref{XXX4.5}-\ref{XXX4.6}).


\begin{Cor} \label{XXX2.10} A regular space $X$ is metrizable
if it has a $\sigma$-locally finite basis of open sets.

\end{Cor}

\begin{pf}  Case 1. $X$ is normal. Suppose $\{U_{s,n}\}_{(s,n)\in S\times N}$
is a basis of open sets in $X$ such that
$\{U_{s,n}\}_{s\in S}$ is locally finite for each $n$.
We may assume that $\{U_{s,n}\}_{s\in S}$ is a cover of $X$ for each $n\in N$.
Notice that each open set $U$ is the union of
countably many of its subsets $F_n$, where $F_n$ is the union
of closures of those $U_{s,n}$ so that $cl(U_{s,n})\subset U$.
Therefore each $U_{s,n}$ is equal to $f_{s,n}^{-1}(0,1]$
for some continuous function $f_{s,n}:X\to [0,1]$.
Notice that $f_n:=\sum\limits_{s\in S}f_{s,n}$
is continuous and maps $X$ to $[0,\infty)$.
Replacing $f_{s,n}$ by $f_{s,n}/f_n$ we may assume
that $f_n\equiv 1$.
Define $g_{s,n}$ as $f_{s,n}\cdot 2^{-n}$
and notice that it induces a partition of unity on $X$
such that $U_{s,n}=g_{s,n}^{-1}(0,1]$ for all $(s,n)\in S\times N$.
By \ref{XXX0.15} the space $X$ is metrizable.
\par To show that $X$ is normal let us observe that,
for any two disjoint, closed subsets $A$ and $B$ of $X$,
there is a countable family $\{U_n\}_{n=1}^\infty$
of open sets in $X$ covering $A$ such that $B\cap cl(U_n)=\emptyset$
for each $n$. Indeed, $U_n$ can be defined as the union
of those $U_{s,n}$ so that $B\cap cl(U_{s,n})=\emptyset$.
Similarly, there is a countable family $\{V_n\}_{n=1}^\infty$
of open sets in $X$ covering $B$ such that $A\cap cl(V_n)=\emptyset$
for each $n$. 
Let $U'_n:=U_n\setminus \bigcup\limits_{k\leq n} cl(V_k)$
and $V'_n:=V_n\setminus \bigcup\limits_{k\leq n} cl(U_k)$.
Notice that $\{U'_n\}_{n=1}^\infty$ is an open cover of $A$,
$\{V'_n\}_{n=1}^\infty$ is an open cover of $B$,
and $U'_n\cap V'_m=\emptyset$ for all $m,n$.
To verify the disjointness of $U'_n$ and $V'_m$
we may assume $n\le m$ without loss of generality.
Now $U'_n\subset U_n$ and $V'_m\subset X\setminus cl(U_n)$,
so those sets are in fact disjoint.
Finally, $U:= \bigcup\limits_{k=1}^\infty U'_k$
and $V:= \bigcup\limits_{k=1}^\infty V'_k$
are two disjoint open subsets of $X$ containing $A$ and $B$,
respectively, which proves that $X$ is normal.
\end{pf}

\section{Extensions of partitions of unity. }

The purpose of this section is to provide a proof of Theorem \ref{XXX0.8}.


\begin{Prop} \label{XXX3.1} Suppose $\cal U=\{U_s\}_{s\in S}$
is an open covering of a space $X$,
$f:X\to [0,\infty)$ is a continuous function, 
$V$ is a neighborhood of $f^{-1}(0,\infty)$ in $X$, and
 $\{f_s\}_{s\in S}$ is a partition of $f|V$ which is
$\cal U|V$-small. 
The extensions $g_s$ of $f_s$ so that $g_s(X-f^{-1}(0,\infty))\subseteq \{0\}$
for each $s\in S$ form a partition of $f$ which is $\cal U$-small.

\end{Prop}

\begin{pf}  It suffices to prove that $g_s$ is continuous
for each $s\in S$. That can be reduced to showing that $g_s^{-1}[0,M)$
is a neighborhood of any $a\in f^{-1}(0)$.
Since $f$ is continuous, there is a neighborhood $U$ of $a$ in $X$
so that $f(U)\subseteq [0,M)$. Notice that $g_s\leq f$
which implies $g_s(U)\subseteq [0,M)$.
\end{pf}

The next result means that global extensions of partitions
exist if one has an extension over a neighborhood.


\begin{Prop} \label{XXX3.2} Suppose $\cal U=\{U_s\}_{s\in S}$
is an open covering of a space $X$, $A$ is a closed subset of $X$, 
$f:X\to [0,\infty)$ is a continuous function, and
 $\{f_s\}_{s\in S}$ is a partition of $f| A$ which is
$\cal U| A$-small. If there is a neighborhood $V$
of $A$ in $X$ such that  $\{f_s\}_{s\in S}$ extends to
a $\cal U| V$-small partition $\{h_s\}_{s\in S}$ of $f| V$,
then $\{f_s\}_{s\in S}$ extends to
a $\cal U$-small partition $\{g_s\}_{s\in S}$ of $f$
if one of the following conditions is satisfied:
\par{a.} $S$ is finite and $X$ is normal,
\par{b.} $X$ is paracompact.
\par\noindent Moreover, if $\{h_s\}_{s\in S}$
is locally finite (respectively, $\{h_s|V-A\}_{s\in S}$
is locally finite), then we may require $\{g_s\}_{s\in S}$
to be locally finite (respectively, $\{g_s|X-A\}_{s\in S}$
to be locally finite).

\end{Prop}

\begin{pf}  We show both cases together. Using \ref{XXX1.12},
pick a locally finite partition of unity
$\{r_s\}_{s\in S}$ on $X$ which is $\cal U$-small. 
Find a neighborhood $W$ of $A$ in $X$ so that
the closure of $W$ is contained in $V$.
Choose a continuous function $u:X\to [0,1]$ so that $u(A)\subseteq \{0\}$
and $u(X-W)\subseteq \{1\}$.
Define $g_s:=(1-u)\cdot h_s+f\cdot u\cdot r_s$ for each $s\in S$.
Notice that $g_s(X-U_s)\subseteq\{0\}$ for each $s\in S$
and $f=\sum\limits_{s\in S}g_s$. Other requirements are easy to see.
\end{pf}

Here is a generalization of \ref{XXX0.4}.


\begin{Lem} \label{XXX3.3} Suppose $X$ is a normal space,
$f:X\to [0,\infty)$ is continuous,
$A$ is a closed subset of $X$, and $\cal U=\{U_s\}_{s\in S}$ is a finite open cover of $X$.
Any $\cal U|A$-small partition $\{f_s\}_{s\in S}$ of $f|A$
extends to
a $\cal U$-small partition $\{g_s\}_{s\in S}$ of $f$.

\end{Lem}

\begin{pf}  
First consider the case of $f$ being positive-valued.
For each $s\in S$ find an extension $h_s:X\to [0,\infty)$
of $f_s$ such that $h_s(X-U_s)\subseteq \{0\}$. 
The continuous function $h:=\sum\limits_{s\in S}h_s:X\to R$
is positive on some neighborhood $W$ of $A$.
Put $p_s:=f\cdot h_s/h$ on $W$ for each $s\in S$.
$\{p_s\}_{s\in S}$ is a partition of $f|W$
and is $\cal U|W$-small. Applying \ref{XXX3.2}
one constructs a partition $\{g_s\}_{s\in S}$ of $f$
on $X$ such that $\{g_s\}_{s\in S}$ is $\cal U$-small.
\par Let $V=f^{-1}(0,\infty)$. By \ref{XXX0.7}, $V$ is normal. Hence,
 $\{f_s|V\cap A\}_{s\in S}$ extends to
a partition $\{p_s\}_{s\in S}$ of $f|V$
which is $\cal U|V$-small. Applying \ref{XXX3.1} one gets that $\{f_s\}_{s\in S}$
extends to
a partition $\{g_s\}_{s\in S}$ of $f$
on $X$ such that $\{g_s\}_{s\in S}$ is $\cal U$-small.
\end{pf}


\begin{Lem} \label{XXX3.4} Suppose $X$ is a paracompact space,
$f:X\to [0,\infty)$ is continuous,
$A$ is a closed subset of $X$, and $\cal U=\{U_s\}_{s\in S}$ is an open cover of $X$.
Any $\cal U|A$-small, locally finite partition $\{f_s\}_{s\in S}$ of $f|A$
extends to
a $\cal U$-small, locally finite partition $\{g_s\}_{s\in S}$ of $f$.

\end{Lem}

\begin{pf}  For each $a\in A$ there is an open subset
$V_a$ of $X$ such that $\{f_s|V_a\cap A\}_{s\in S}$
is finite. By adding $X-A$ one creates an open cover
$\cal V=\{V_p\}_{p\in P}$ of $X$ with the property
that $\{f_s|V_p\cap A\}_{s\in S}$ is finite
for each $p\in P$. 
\par
Choose a locally finite partition of unity
$\{h_p\}_{p\in P}$ on $X$ which is $\cal V$-small (use \ref{XXX1.12}).
For each finite subset $T$ of $P$ define $B_T$
as $\{x\in X\mid h_p(x)>0\ \implies \ p\in T\}$.
Notice that $B_T$ is closed (if $h_p(x)>0$ for some $p\notin T$,
then the neighborhood $\{y\in X\mid h_p(y)>0\}$ is contained
in $X\setminus B_T$) and $B_F\subseteq B_T$ if $F\subseteq T$.
We plan to create, by induction on the size of $T$,
a finite partition $\{f^T_s\}_{s\in S}$
of $f|B_T$ which is $\cal U|B_T$-small, extends
$\{f_s|A\cap B_T\}_{s\in S}$,
and $\{f^T_s\}_{s\in S}$ extends $\{f^F_s\}_{s\in S}$
for $F\subseteq T$.
\par
Since $B_T\subseteq \bigcup\limits_{p\in T}V_p$,
$\{f_s|A\cap B_T\}_{s\in S}$ is finite for each
finite $T\subseteq P$.
Notice that $B_p$, $p\in P$, are mutually disjoint,
so using \ref{XXX3.3} one can create $\{f^p_s\}_{s\in S}$
for each $p\in P$. Once $\{f^F_s\}_{s\in S}$
exist for all $F$ containing less that $n$ elements,
$\{f^T_s\}_{s\in S}$ (for any $T$ containing $n$ elements)
can be constructed by pasting $\{f^F_s\}_{s\in S}$
for all $F\subset T$ with $\{f_s|A\cap B_T\}_{s\in S}$,
and then extending over $B_T$ using \ref{XXX3.3}.
\end{pf}


\begin{Lem} \label{XXX3.5} Suppose $\{f_s\}_{s\in S}$ is a partition of 
a continuous function $f:X\to (0,\infty)$. 
There exists a partition $\{f_n\}_{n=1}^\infty$
of $f$ and there exist locally finite partitions $\{f_s^n\}_{s\in S}$ of $f_n$
such that $\{f_s^n\}_{n=1}^\infty$ is a partition of $f_s$ for each $s\in S$.

\end{Lem}

\begin{pf}  Let $g=\sup\{f_s\mid s\in S\}$.
Put $f^1_s=\max(0,f_s-g/2)$ and $h^1_s=f_s-f^1_s$.
Notice that $\sup\{f^1_s\mid s\in S\}=g/2=\sup\{h^1_s\mid s\in S\}$
and $\{f_s^1\}_{s\in S}$ is a locally finite partition of some $f_1$
(see \ref{XXX1.8} and the proof of \ref{XXX1.12}).
Apply the same step to $\{h^1_s\}_{s\in S}$ and extract 
a locally finite partition $\{f^2_s\}_{s\in S}$ of $f_2$
such that $\sup\{f^2_s\mid s\in S\}=g/4=\sup\{h^1_s-f_s^2\mid s\in S\}$.
Continuing inductively one expresses each $f_s$ as
$\sum\limits_{n=1}^\infty f_s^n$ so that $\{f^n_s\}_{s\in S}$ is
a locally finite partition of $f_n$.
Clearly, $\{f_n\}_{n=1}^\infty$ is a partition of $f$.
\end{pf}


\begin{Lem} \label{XXX3.6} Suppose $X$ is a paracompact space,
$f:X\to [0,\infty)$ is continuous,
$A$ is a closed subset of $X$, and $\cal U=\{U_s\}_{s\in S}$ is an open cover of $X$.
Any $\cal U|A$-small partition $\{f_s\}_{s\in S}$ of $f|A$ extends to
a $\cal U$-small partition $\{g_s\}_{s\in S}$ of $f$.

\end{Lem}

\begin{pf}  
First consider the case of $f>0$.
This case can be farther reduced to that of $f\equiv 1$ by switching
to $\{f_s/f\}_{s\in S}$, extending it over $X$, and multiplying
the extension by $f$.
By \ref{XXX3.5} there exists a partition of unity $\{f_n\}_{n=1}^\infty$
on $A$ and there exist locally finite partitions $\{f_s^n\}_{s\in S}$ of $f_n$
such that $\{f_s^n\}_{n=1}^\infty$ is a partition of $f_s$ for each $s\in S$.
Extend $\{f_n\}_{n=1}^\infty$ to a partition of unity
$\{g_n\}_{n=1}^\infty$ on $X$ (such an extension exists by \ref{XXX0.2},
see also 6.5 of \cite{D$_1$}).
Notice that $\{f_s^n\}_{s\in S}$ is $\cal U|A$-small and,
by \ref{XXX3.4}, there is a locally finite partition $\{g_s^n\}_{s\in S}$ 
of $g_n$ on $X$ which is $\cal U$-small.
By \ref{XXX1.5}, $g_s=\sum\limits_{n=1}^\infty g^n_s$ is continuous and, clearly,
$\{g_s\}_{s\in S}$ is a $\cal U$-small partition of unity extending
$\{f_s\}_{s\in S}$.
\par Let $V=f^{-1}(0,\infty)$. $V$ is paracompact by \ref{XXX0.7},
so $\{f_s|V\cap A\}_{s\in S}$ extends to
a partition $\{p_s\}_{s\in S}$ of $f|V$
which is $\cal U|V$-small. Applying \ref{XXX3.1} one gets that $\{f_s\}_{s\in S}$
extends to
a partition $\{g_s\}_{s\in S}$ of $f$
on $X$ such that $\{g_s\}_{s\in S}$ is $\cal U$-small.
\end{pf}


\begin{Lem} \label{XXX3.7} Suppose $X$ is a paracompact space,
$A$ is a closed G$_\delta$-subset of $X$, and $\cal U=\{U_s\}_{s\in S}$ is an open cover of $X$.
Any partition $\{f_s\}_{s\in S}$ of unity on $A$
which is $\cal U|A$-small extends to
a partition $\{g_s\}_{s\in S}$ of unity
on $X$ such that $\{g_s\}_{s\in S}$ is $\cal U$-small
and $\{g_s|X-A\}_{s\in S}$ is locally finite.

\end{Lem}

\begin{pf}  Pick $u:X\to [0,1]$ so that $A=u^{-1}(0)$.
Extend $\{f_s\}_{s\in S}$ to a partition of unity $\{h_s\}_{s\in S}$
on $X$ which is
$\cal U$-small (see \ref{XXX3.6}).
Define $V_s$ as $\{x\in U_s\mid h_s(x)>u(x)\}$
and notice that $V=\bigcup\limits_{s\in S} V_s$
contains $A$. Also, $\{f_s\}_{s\in S}$ is $\cal V|A$-small,
where $\cal V=\{V_s\}_{s\in S}$. 
By \ref{XXX3.6} there is an extension $\{h'_s\}_{s\in S}$
of $\{f_s\}_{s\in S}$ over $V$ which is $\cal V$-small.
Notice that $\{h'_s|V-A\}_{s\in S}$ is locally finite
(see \ref{XXX1.5} and \ref{XXX1.8}).
By \ref{XXX3.2} one can construct a partition of unity $\{g_s\}_{s\in S}$
on $X$ which extends $\{f_s\}_{s\in S}$ so that $\{g_s\}_{s\in S}$ is $\cal U$-small
and $\{g_s|X-A\}_{s\in S}$ is locally finite.
\end{pf}

\section{Integrals and derivatives of partitions of unity. }


\begin{Def} \label{XXX4.1} Suppose $S\ne\emptyset$ is a set, $X$ is a space, and $\{f'_T\}_{T\subseteq S}$
is a partition of unity on $X$ indexed by all finite subsets $T\ne\emptyset$ of $S$.
The {\bf integral} of $\{f'_T\}_{T\subseteq S}$ is the partition
of unity $\{f_s\}_{s\in S}$ defined as follows:
$f_s(x)$ is the sum of all $f'_T(x)/|T|$, where $s\in T$
and $|T|$ is the number of elements of $T$. 

\end{Def} 

Notice that each $f_s$ is continuous by \ref{XXX1.5}. Also, it is clear that
$\sum\limits_{s\in S}f_s=\sum\limits_{T\subseteq S}f'_T$,
so $\{f_s\}_{s\in S}$ is indeed a partition of unity on $X$.


\begin{Thm} \label{XXX4.2} For any partition of unity $\{f_s\}_{s\in S}$ on a space $X$
there is a unique
partition of unity $\{f'_T\}_{T\subseteq S}$ on $X$ (called the {\bf derivative}
of $\{f_s\}_{s\in S}$)
satisfying the 
following properties:
\par{1.} $\{f_s\}_{s\in S}$ is the integral of $\{f'_T\}_{T\subseteq S}$.
\par{2.} $f'_T(x)\ne 0$ and $f'_F(x)\ne 0$ for some $x\in X$
implies $T\subseteq F$ or $F\subseteq T$. 
\par\noindent It is given by
$f'_T=|T|\cdot \max(0,g_T)$, where
$g_T=\min\{f_t\mid t\in T\}-\sup\{f_t\mid t\in S-T\}$. 

\end{Thm}

\begin{pf}   Suppose $\{f'_T\}_{T\subseteq S}$ 
has properties 1 and 2. Notice that, for each $x\in X$,
$\{T\subset S:f'_T(x)>0\}$ is countable.
Given $x\in X$ we can, by using 2, find a (possibly finite) strictly increasing
sequence
$T(1)\subset T(2)\subset \ldots$
of finite subsets of $S$ such that $f'_T(x)>0$ if and only if
$T$ equals $T(i)$ for some $i$.
Let $M$ be the supremum of all $i$ such that
$T(i)$ exists (it is possible for $M$ to be $\infty$).
For integers $i\leq M$ let $v_i=\sum\limits_{k=i}^M f'_{T(k)}(x)/|T(k)|$.
Notice that it is a strictly decreasing sequence of positive numbers.
The meaning of those numbers is as follows:
$v_1=f_s(x)$ for $s\in T(1)$, and
$v_{i+1}=f_s(x)$ for $x\in T(i+1)\setminus T(i)$ if $i\ge 1$.
Notice that $f_s(x)=0$ for $s\notin \bigcup\limits_{i=1}^M T(i)$.
Define $g_T(x)=\min\{f_t(x)\mid t\in T\}-\sup\{f_t(x)\mid t\in S-T\}$
and $h_T(x)=|T|\cdot \max(0,g_T(x))$ for all finite subsets $T$ of $S$.
To prove uniqueness of the derivative we need to show
$f'_T(x)=h_T(x)$ for all $T$.
There are two possible cases:
\par A. $T\cap T(i)=\emptyset$ for all $i\leq M$.
\par B. $T\cap T(i)\ne\emptyset$ for some $i\leq M$.
\par In Case A one has $f_s(x)=0$ for all $s\in T$ implying
$g_T(x)\leq 0$, and $f'_T(x)=0$. Thus, $f'_T(x)=h_T(x)=0$
in that case.
\par In Case B, if $T\setminus \bigcup\limits_{i=1}^M T(i)$
contains some $s$, then $f'_T(x)=0$ and $g_T(x)\leq f_s(x)=0$,
so again $f'_T(x)=h_T(x)=0$.
Therefore, one may pick the smallest $i$ such that $T\subseteq T(i)$.
If there is $s\in T(i)\setminus T$, then
$f'_T(x)=0$ (as $T$ does not equal to any of $T(j)$)
and $f_s(x)\ge f_t(x)$ for all $t\in T$, implying
$g_T(x)\leq 0$. Again, $f'_T(x)=h_T(x)=0$.
Thus, $T=T(i)$. Pick $s\in T(i)\setminus T(i-1)$
($s\in T(1)$ if $i=1$). Now $f_s(x)=v_i$ is the smallest value
of all $f_t(x)$, $t\in T$. Also, $v_{i+1}$ ($0$, if $i=M$)
is the largest value of $f_t(x)$, $t\in S\setminus T$.
That means $h_T(x)=|T|(v_i-v_{i+1})=f'_T(x)$
and we are done with the proof of uniqueness.
\par To prove existence of the derivative, put $f'_T=|T|\cdot \max(0,g_T)$, where
$g_T=\min\{f_t\mid t\in T\}-\sup\{f_t\mid t\in S-T\}$.
By \ref{XXX1.9}-\ref{XXX1.10}, $g_T$ is continuous which implies continuity of $f'_T$.
Suppose $g_T(x)>0$ and $g_F(x)>0$ for some $x\in X$.
If $s\in T-F$, then $f_s(x)>f_t(x)$ for all $t\in S-T$.
Similarly, $t\in F-T$ implies $f_t(x)>f_s(x)$, which means such $t$
does not exist and $F\subseteq T$.
It remains to show that the sum of all $\max(0,g_T)$ with $s\in T$ is $f_s$ if $s$ is fixed.
Indeed, we can pick elements $s(i)\in S$, $i\ge 1$ or $n\ge i\ge 1$
for some $n$,
such that $f_{s(i)}(x)\ge f_{s(i+1)}(x)>0$ for all $i$ and $f_t(x)=0$ for $t$ not equal
to any of $s(i)$. 
\par If $s\ne s(i)$ for all $i$, then $f_s(x)=0$
implying $g_T(x)=0$ for all $T$ containing $s$.
Consequently, in that case, $f_s(x)$ is the sum of 
all $\max(0,g_T)$ with $s\in T$.
\par Suppose $s=s(i)$ for some $i$. Now, $g_T(x)>0$
and $s\in T$ can happen only if $T=T(k)=\{s(1),\ldots,s(k)\}$ for
some $k\ge i$.
Notice that $g_{T(k)}(x)=f_{s(k)}(x)-f_{s(k+1)}(x)$ in that case.
Now $\sum\limits_{k=i}^\infty g_{T(k)}(x)$ becomes a telescopic sequence
 $(f_{s(i)}(x)-f_{s(i+1)}(x))+(f_{s(i+1)}(x)-f_{s(i+2)}(x))+\ldots$
which adds up to $f_s(x)$.  
\par It remains to show that $\{f'_T\}_{T\subseteq S}$ is actually a
partition of unity on $X$. Notice that the sum 
$\sum\limits_{T\subseteq S}f'_T(x)$
can be expressed as $\sum\limits_{s\in S}\sum\limits_{s\in T\subseteq S}f'_T(x)/|T|$
which is $\sum\limits_{s\in S}f_s(x)=1$.
\end{pf}

Let us show how to use derivatives of partitions of unity
to create star refinements of open covers, a basic operation in general topology.


\begin{Prop} \label{XXX4.3} Suppose $\{f_s\}_{s\in S}$ is a partition of unity  
on $X$ and $\{f'_T\}_{T\subseteq S}$ is its
derivative. 
If open covers $\cal U=\{U_T\}_{T\subseteq S}$
and $\cal V=\{V_s\}_{s\in S}$ of $X$ are defined by
$U_T=(f'_T)^{-1}(0,1]$ and $V_s=f_s^{-1}(0,1]$,
then stars of $\cal U$ at points of $X$ refine $\cal V$. 

\end{Prop}

\begin{pf}   Given $x\in X$ let us pick $s\in S$ so that
$f_s(x)=\sup\{f_t(x)\mid t\in S\}$ and let $F=\{t\in S\mid f_t(x)=f_s(x)\}$.
Since $f'_T=|T|\cdot \max(0,g_T)$, where
$g_T=\min\{f_t\mid t\in T\}-\sup\{f_t\mid t\in S-T\}$,
$f'_F(x)>0$ and $f'_T(x)=0$ for every proper subset $T$ of $F$.
Suppose $x\in U_T$ for some $T\subseteq S$. That means $f'_T(x)>0$
and $T$ must contain $F$. Hence, $s\in T$
which implies $f_s(y)\ge f'_T(y)/|T|>0$ for all $y\in U_T$.
That proves that the star of $\cal U$ at $x$ is contained in $V_s$. 
\end{pf}

The following allows to apply the calculus
of partitions of unity in dimension theory.


\begin{Prop} \label{XXX4.4} Let $\{f_s\}_{s\in S}$ be a partition of unity on a space $X$
and let $\{f'_T\}_{T\subseteq S}$ be its derivative.
The order of $\{f_s\}_{s\in S}$ is at most $n$ if and only if
$f'_T\equiv 0$ for all $T\subseteq S$ containing at least $(n+2)$ elements. 

\end{Prop}

\begin{pf}   Suppose the order of $\{f_s\}_{s\in S}$ is at most $n$
and suppose $T$ is a subset of $S$ containing at least $n+2$ elements.
Since $f'_T=|T|\cdot \max(0,g_T)$, where
$g_T=\min\{f_t\mid t\in T\}-\sup\{f_t\mid t\in S-T\}$,
and since
for any $x\in X$ there must be at least one $s\in S$ with $f_s(x)=0$,
we get $f'_T\equiv 0$.
\par Suppose $f'_T\equiv 0$ for all $T\subseteq S$ containing at least $n+2$ elements.
Suppose there is $F\subseteq S$ containing at least $n+2$ elements such that for some
$x\in X$ all values $f_s(x)$, $s\in F$, are positive.
Let $a=\min\{f_s(x)\mid s\in F\}$. Enlarge $F$, if necessary,
to include all $s\in S$ such that $f_s(x)\ge a$.
Now, $f'_F(x)=|F|\cdot(\min\{f_s(x)\mid s\in F\}-\sup\{f_s(x)\mid s\in S-F\})>0$,
a contradiction. 
\end{pf}

Another application of derivatives of partitions of unity yields
the second part of the classical metrizability criterion.


\begin{Cor}[A.H.Stone \cite{En$_1$}] \label{XXX4.5}  Each open covering of a
 metrizable space $X$
has a $\sigma$-discrete refinement.

\end{Cor}

\begin{pf}   It suffices to show
that, for any partition of unity $f=\{f_s\}_{s\in S}$ on $X$,
the cover $\{U_s\}_{s\in S}$ of $X$,
$U_s:=f_s^{-1}(0,1]$, has a $\sigma$-discrete refinement.
Let $\{f'_T\}_{T\subseteq S}$ be the derivative of $f$.
By considering only $T$ of a given size $n$
one gets that the sets $U_T:=\{x\in X| f'_T(x)>0\}$
are mutually disjoint. Therefore, by \ref{XXX1.4} and \ref{XXX1.8},
the family $\cal U_{m,n}$ consisting of sets
$U_{T,m}:=\{x\in X| f'_T(x)>1/m\}$, where $|T|=n$,
is locally finite and the closures of its elements are mutually disjoint.
That means precisely that $\cal U_{m,n}$ is discrete.
\end{pf}


\begin{Cor}[Bing-Nagata-Smirnov \cite{En$_1$}] \label{XXX4.6}  A regular
space $X$ is metrizable
if and only if it
has a $\sigma$-discrete basis of open sets.

\end{Cor}

\section{Dimension and partitions of unity. }

\par

Using partitions of unity one can introduce the covering
dimension of normal spaces via finite partitions of unity,
and for paracompact spaces via arbitrary partitions of unity.
We chose arbitrary partitions of unity in order to illustrate
how one applies the calculus of partitions of unity.
One can show (using techniques of \cite{D$_3$}) that both ways
yield the same result for paracompact spaces.


\begin{Def} \label{XXX5.1} Let $\cal U$ be an open cover of a space $X$.
The order $ord(\cal U)$ of $\cal U$ is the smallest integer $n$ with the property
that any family $U_1,\ldots,U_{n+2}$ of different elements of $\cal U$
has empty intersection. 

\end{Def}


\begin{Rem} \label{XXX5.2}  Notice that if $\{f_s\}_{s\in S}$ is a partition of unity on $X$,
then its order is the same as that of the open covering $\{U_s\}_{s\in S}$,
where $U_s=f_s^{-1}(0,1]$ for each $s\in S$. 

\end{Rem} 

\par


\begin{Lem} \label{XXX5.3} Let $n\ge 0$. If $X$ is a paracompact space, then the following conditions
are equivalent:
\par{1.} Any open cover $\cal U$ of $X$ has an open refinement $\cal V$
such that $ord(\cal V)\leq n$.
\par {2.} For any open cover $\{U_s\}_{s\in S}$ of $X$ there is
a partition of unity $\{f_s\}_{s\in S}$ of order at most $n$ such that
$f_s(X-U_s)\subseteq \{0\}$ for each $s\in S$.
\par {3.} Any partition of unity $\{f_s\}_{s\in S}$ on $X$
is approximable by partitions of unity of order at most $n$. 

\end{Lem}

\begin{pf}   1)$\implies$2). Let $\{U_s\}_{s\in S}$ be an open cover
of $X$. Pick a refinement $\{V_t\}_{t\in T}$ of $\{U_s\}_{s\in S}$ whose order
is at most $n$. Let $\{g_t\}_{t\in T}$ be a partition of unity
on $X$ such that $g_t(X-V_t)\subseteq \{0\}$ for each $t\in T$.
Notice that the order of $\{g_t\}_{t\in T}$ does not exceed $n$.
Partition $T$ into disjoint subsets $T_s$, $s\in S$,
such that $g_t(X-U_s)\subseteq \{0\}$ for all $t\in T_s$.
If $T_s=\emptyset$ we put $f_s=0$, otherwise $f_s=\sum\limits_{t\in T_s}g_t$.
Notice that $\{f_s\}_{s\in S}$ is the required partition of unity
on $X$.
\par 2)$\implies$3). Given a partition of unity $\{f_s\}_{s\in S}$
on $X$, put $V_s=f_s^{-1}(0,1]$ and find a partition
of unity $\{g_s\}_{s\in S}$ on $X$ of order at most $n$
such that $g_s(X-V_s)\subseteq \{0\}$ for each $s\in S$.
Clearly, $\{g_s\}_{s\in S}$ approximates $\{f_s\}_{s\in S}$.
\par 3)$\implies$1). Given open cover $\cal U=\{U_s\}_{s\in S}$ of $X$
let us pick a partition of unity $\{f_s\}_{s\in S}$ such that
$f_s(X-U_s)\subseteq \{0\}$ for each $s\in S$.
Approximate $\{f_s\}_{s\in S}$ by $\{g_s\}_{s\in S}$ which is of
order at most $n$.
Put $V_s=g_s^{-1}(0,1]$ and notice that $\{V_s\}_{s\in S}$ refines
$\{U_s\}_{s\in S}$ and its order does not exceeed $n$. 
\end{pf}


\begin{Def} \label{XXX5.4} Let $X$ be a paracompact space and $n\ge -1$.
$\dim(X)=-1$ means that $X$ is empty. Suppose $n\ge 0$.
We say that $X$ is at most $n$-dimensional (notation: $\dim(X)\leq n$)
if one of the conditions 1, 2, or 3 of \ref{XXX5.3} is satisfied.
We say that $X$ is $n$-dimensional (notation: $\dim(X)=n$) if $\dim(X)\leq n$
and $\dim(X)\leq n-1$ does not hold. 

\end{Def} 

\par


\begin{Cor}[\cite{En$_2$}, 3.1.3] \label{XXX5.5}  Suppose $X$ is a paracompact space.
If $\dim(X)\leq n$ and $A$ is a closed subset
of $X$, then $\dim(A)\leq n$. 

\end{Cor}

\begin{pf}   Suppose $\{U_s\}_{s\in S}$ is an open cover of $A$.
Define $V_s=(X-A)\cup U_s$ for $s\in S$ and notice that $V_s$ is an open subset of
$X$. Since $\{V_s\}_{s\in S}$ is an open cover of $X$, there is
a partition of unity $\{f_s\}_{s\in S}$ on $X$ such that $f_s(X-V_s)\subseteq \{0\}$
for each $s\in S$ and its order does not exceed $n$. Since $X-V_s=A-U_s$, the partition of unity
$\{g_s\}_{s\in S}=\{f_s|A\}_{s\in S}$ on $A$ satisfies $g_s(A-U_s)\subseteq \{0\}$
and its order does not exceed $n$. By \ref{XXX5.3}, $\dim(A)\leq n$. 
\end{pf}
 

\begin{Proof} \label{XXXpfE} Proof of \ref{XXX0.11}.
\end{Proof}
\begin{pf}
  By \ref{XXX3.6} there is a partition of unity 
 $\{h_s\}_{s\in S}$ 
on $X$
such that $h_s(X-U_s)\subseteq \{0\}$
and $h_s|A=f_s$ for each $s\in S$. 
Let $\{h'_T\}_{T\subseteq S}$ be its derivative.
Consider the sum $h$ of all $h'_T$ such that $T$ contains
at most $(n+1)$ elements. By \ref{XXX4.4}, $h|A$ is equal $1$.
Let $W=\{x\in X\mid h(x)\ne 0\}$. 
Define $i'_T$ on $W$ as $h'_T/h$ if $T$ contains at most $(n+1)$ elements.
Otherwise put $i'_T=0$.
Integrate it to $\{i_s\}_{s\in S}$ and notice that
$\{i'_T\}_{T\subseteq S}$ is its derivative.
By \ref{XXX4.4}, the order of $\{i_s\}_{s\in S}$ is at most $n$.
Also, it is clear that $i_s(W-U_s)\subseteq \{0\}$ for each $s\in S$.
Pick an open subset $V$ containing $A$ whose closure is contained in $W$.
Pick an open subset $U$ containing $A$ whose closure is contained in $V$.
If $\dim(X)>n$, then we put $B=cl(V)$ and we extend $\{i_s|B\}_{s\in S}$
over $X$ to obtain $\{g_s\}_{s\in S}$ satisfying conditions a)-c).
If $\dim(X)\leq n$, we may extend $\{i_s|cl(V)\}_{s\in S}$
over $X$ so that $i_s(X-U_s)\subseteq \{0\}$ for each $s\in S$
and
we choose a partition of unity  $\{j_s\}_{s\in S}$ 
on $X$ of order at most $n$ which approximates $\{i_s\}_{s\in S}$ (see \ref{XXX5.5} and \ref{XXX5.3}).
Let $a:X\to [0,1]$ be a continuous function such that $a(cl(U))\subseteq \{1\}$ and $a(X-V)\subseteq \{0\}$.
Define $g_s(x)$ as $a(x)\cdot i_s(x)+(1-a(x))\cdot j_s(x)$. 
Notice that $g_s$ is continuous for each $s\in S$,
the sum $\sum\limits_{s\in S}g_s$ equal $1$,
and $\{g_s\}_{s\in S}$ satisfies conditions a)-b).
Pick $x\in X$ and suppose that $T=\{s\in S\mid g_s(x)\ne 0\}$
contains more that $n+1$ elements. This can happen only if $0<a(x)<1$,
in particular $x\in V$.
Since $j_s(x)\ne 0$ implies $i_s(x)\ne 0$, we get $i_s(x)>0$ for all $s\in T$
contradicting the fact that the order of $\{i_s\}_{s\in S}$ on $V$
is at most $n$. Thus, conditions c) and d) hold. 
\end{pf}


\begin{Thm}[\cite{En$_2$}, 3.1.8] \label{XXX5.7}  Suppose $X$ is a paracompact space and $X=\bigcup\limits_{k=1}^\infty X_k$,
where $X_k$ is closed in $X$ for each $k$.
If $\dim(X_k)\leq n$ for each $k$, then $\dim(X)\leq n$ . 

\end{Thm}

\begin{pf}   Suppose $\{f_s\}_{s\in S}$ is a partition of
unity on $X$. Approximate $\{f_s|X_1\}_{s\in S}$
by a partition of unity of order at most $n$ and extend it over
$X$ so that the resulting partition of unity $\{g_{1,s}\}_{s\in S}$
approximates $\{f_s\}_{s\in S}$ and the order of $\{g_{1,s}|B_1\}_{s\in S}$
is at most $n$ for some closed neighborhood $B_1$ of $X_1$ in $X$.
Suppose that, for some $k\ge 1$, there is a partition of unity
$\{g_{k,s}\}_{s\in S}$ on $X$ which approximates $\{f_{s}\}_{s\in S}$
and, for some closed neighborhood $B_k$ of $\bigcup\limits_{i=1}^k X_k$,
the order of $\{g_{k,s}|B_k\}_{s\in S}$ is at most $n$.
Put $A=B_k\cap X_{k+1}$. Since $\dim (A)\leq n$ (see \ref{XXX4.5}),
$\{g_{k,s}|A\}_{s\in S}$ can be extended over $X_{k+1}$
to approximate $\{f_{s}|X_{k+1}\}_{s\in S}$ and preserving
the order at the same time. Pasting the extension with $\{g_{k,s}|B_k\}_{s\in S}$
and then extending over $X$ using \ref{XXX0.11} gives an approximation
$\{g_{k+1,s}\}_{s\in S}$ of $\{f_{s}\}_{s\in S}$
so that its order on some closed neighborhood $B_{k+1}$ of $B_k\cup X_{k+1}$
is at most $n$.
The direct limit of all $\{g_{k,s}\}_{s\in S}$ as $k\to\infty$
gives an approximation $\{g_{s}\}_{s\in S}$ of
$\{f_{s}\}_{s\in S}$ whose order is at most $n$.
By \ref{XXX5.3}, $\dim(X)\leq n$. 
\end{pf}

Part of the meaning of \ref{XXX0.11} is that partitions of unity on closed subsets
of paracompact spaces can be extended over a neighborhood while
preserving order. The next result deals with approximate extensions.


\begin{Prop} \label{XXX5.8} Suppose $A$ is a subset of a metrizable space
$X$ and $\{f_{s}\}_{s\in S}$ is a partition of unity on $A$.
There
is a neighborhood $U$ of $A$ in $X$ and a locally finite partition of unity
$\{g_{s}\}_{s\in S}$ on $U$ so that 
$\{g_{s}|A\}_{s\in S}$
approximates $\{f_{s}\}_{s\in S}$.
Moreover, if order of $\{f_{s}\}_{s\in S}$ is at most $n$, then 
we may require
$\{g_{s}\}_{s\in S}$ to be of order at most $n$. 

\end{Prop}

\begin{pf}   Given an open set $U$ of $A$
define $e(U)$ as $\{x\in X\mid dist(x,A)<dist(x,X\setminus U)\}$.
If $x\in A\cap e(U)$, then $0=dist(x,A)<dist(x,X\setminus U)$,
i.e. $x\in U$. Conversely, $x\in A\cap U$
implies $0=dist(x,A)<dist(x,X\setminus U)$, i.e. $x\in e(U)$.
Notice that
$e(V\cap W)=e(V)\cap e(W)$ for any two open subsets $V$ and $W$
of $A$. Indeed, it follows from the equality $dist(x,X\setminus V\cap W)=
\min(dist(x,X\setminus V),dist(x,X\setminus W))$.
\par
Define $U_s:=f_{s}^{-1}(0,1]$ for $s\in S$.
Let $V_s:=e(U_s)$ for each $s\in S$. 
Put $U=\bigcup\limits_{s\in S}V_s$.
Since $e(\bigcap\limits_{s\in T}U_s)=\bigcap\limits_{s\in T}e(U_s)$ for any 
finite subset $T$ of $S$, the order of $\{V_{s}\}_{s\in S}$ is at most 
that of $\{f_{s}\}_{s\in S}$.
Choose a locally finite partition of unity $\{g_{s}\}_{s\in S}$ on $U$
such that $g_{s}(U-V_s)\subseteq \{0\}$ for $s\in S$.
Notice that the order of $\{g_{s}\}_{s\in S}$ is at most that of $\{f_{s}\}_{s\in S}$
and $\{g_{s}|A\}_{s\in S}$ approximates $\{f_{s}\}_{s\in S}$. 
\end{pf}


\begin{Cor}[\cite{En$_2$}, 3.1.23] \label{XXX5.9}  Suppose $A$ is a subset of a space $X$.
If $X$ is metrizable, then
$\dim(A)\leq \dim(X)$. 

\end{Cor}

\begin{pf}  Let $\dim(X)=n$. Given a partition of unity $\{f_{s}\}_{s\in S}$
on $A$ we may find an open neighborhood $U$ of $A$ in $X$
and a partition of unity $\{g_{s}\}_{s\in S}$ on $U$
such that $\{g_{s}|A\}_{s\in S}$ approximates $\{f_{s}\}_{s\in S}$.
Since $U$ is a $F_\sigma$-set in $X$, $\dim(U)\leq n$ by \ref{XXX5.7}
and $\{g_{s}\}_{s\in S}$ is approximable by
$\{h_{s}\}_{s\in S}$ of order at most $n$ (see \ref{XXX5.3}).
Notice that $\{h_{s}|A\}_{s\in S}$ approximates
$\{f_{s}\}_{s\in S}$ and its order is at most $n$.
By \ref{XXX5.3}, $\dim(A)\leq n$. 
\end{pf}


\begin{Thm}[\cite{En$_2$}, 4.1.18] \label{XXX5.10} Suppose $A$ and $B$ are subsets of a space $X$.
If $X$ is metrizable, then
$\dim(A\cup B)\leq \dim(A)+\dim(B)+1$. 

\end{Thm}

\begin{pf}  Let $\dim(A)=m$ and $\dim(B)=n$.
Suppose $\{f_{s}\}_{s\in S}$ is a partition of unity on $X$.
By \ref{XXX5.8} and \ref{XXX5.3} we may find open neighborhoods $U$ of $A$
and $V$ of $B$ such that there exist partitions of unity
$\{g_{s}\}_{s\in S}$ on $U$ of order at most $m$,
and $\{h_{s}\}_{s\in S}$ on $V$ of order at most $n$
such that $\{g_{s}\}_{s\in S}$ approximates $\{f_{s}|U\}_{s\in S}$
and $\{h_{s}\}_{s\in S}$ approximates $\{f_{s}|V\}_{s\in S}$.
Choose a continuous function $a:X\to [0,1]$ such that $a(X-U)\subseteq \{0\}$
and $a(X-V)\subseteq \{1\}$.
Define $p_s(x)=a(x)\cdot g_s(x)+(1-a(x))\cdot h_s(x)$
for $x\in X$. Notice that $\{p_{s}\}_{s\in S}$
is a partition of unity approximating $\{f_{s}\}_{s\in S}$
and whose order is at most $m+n+1$.  
\end{pf}

\section{Simplicial complexes. }

There are two ways of introducing simplicial complexes.
One is abstract and follows the way nerves of open covers
are introduced (see \cite{D-S} or \cite{M-S}).


\begin{Def} \label{XXX6.1} Given a partition of unity $f=\{f_s\}_{s\in S}$
on a space $X$ its {\bf nerve} $\cal N(f)$ is defined
as the set of all finite subsets $T$ of $S$
with the property that there is $x\in X$
with $f_s(x)>0$ for all $s\in T$. 

\end{Def}

Alternatively, the nerve can be defined using derivatives
of partitions of unity and this way is more fruitful.


\begin{Prop} \label{XXX6.2} Suppose $f=\{f_s\}_{s\in S}$ is  a partition of unity
on a space $X$. $T\in \cal N(f)$ if and only if there is
a finite subset $F$ of $S$ containing $T$ such that $f'_F\neq 0$. 

\end{Prop}

\begin{pf}  If $f'_F(x)>0$ for some $F\supset T$, then $f_s(x)\ge f'_F(x)/|F|>0$
for each $s\in T$ which proves $T\in \cal N(f)$. 
Conversely, if there is a point $x\in X$ such that $f_s(x)>0$ for all $s\in T$,
then we put $F:=\{r\in S\mid f_r(x)\geq \min\limits_{s\in T}f_s(x)\}$
and notice that $f'_F(x)=|F|\cdot \max(0,\min\limits_{s\in F}f_s(x)-
\sup\limits_{s\in S-F}f_s(x))>0$.
\end{pf}


\begin{Cor} \label{XXX6.3} Suppose $f=\{f_s\}_{s\in S}$ is  a partition of unity
on a space $X$ and $A$ is a subset of $X$.
There is a neighborhood $U$ of $A$ in $X$ and
  a partition of unity $g=\{g_s\}_{s\in S}$ on $U$
such that the following conditions hold:
\par{a.} $g|A=f|A$.
\par{b.} $g$ approximates $f|U$.
\par{c.} The nerve of $g$ equals the nerve of $f|A$. 
\par{d.} $U$ is the set of all points $x$ 
such that $f'_T(x)>0$ and $f'_T|A\ne 0$ for some finite $T\subset S$.

\end{Cor}

\begin{pf}  Consider the derivative $\{f'_T\}_{T\subseteq S}$ of $f=\{f_s\}_{s\in S}$.
Put $h'_T\equiv 0$ if $f'_T|A\equiv 0$ and $h'_T=f'_T$
otherwise.
Notice that $h=\sum\limits_{T\subseteq S}h'_T$ is continuous
and equals 1 on $A$.
Let $U:=\{x\in X\mid h(x)>0\}$.
Define $g'_T$ as $h'_T/h$. $\{g'_T\}_{T\subseteq S}$
is a partition of unity on $U$. Let $\{g_s\}_{s\in S}$
be its integral. Since $g'_T(x)>0$ and $g'_F(x)>0$ implies
$f'_T(x)>0$ and $f'_F(x)>0$, one gets $F\subset T$ or $T\subset F$.
That means $\{g'_T\}_{T\subseteq S}$
is the derivative of $\{g_s\}_{s\in S}$.
Since $g'_T|A=f'_T|A$ for all finite subsets $T$ of $S$,
$g_s|A=f_s|A$ for all $s\in S$.
\par Suppose $g_s(x)>0$ for some $x\in U$ and $s\in S$.
There is a finite $T\subset S$ containing $s$ such that
$g'_T(x)>0$. Therefore $h'_T(x)>0$ which implies
$h'_T(x)=f'_T(x)$ and $f_s(x)>0$, i.e. $g$ approximates $f|U$.
\par
Suppose $F$ is a finite subset of $S$
containing $T$ and $g'_F\ne 0$. Therefore $h'_F\ne 0$
which means $h'_F=f'_F$ and $f'_F|A\ne 0$. By \ref{XXX6.2}
the nerve of $g$ equals the nerve of $f|A$.
\end{pf}

The second way of introducing simplicial complexes
is much more geometric (see \cite{M-S} for details).
Namely, a simplicial complex is a family $K$
of geometric simplices $\Delta$ with the property
that every face of $\Delta$ belongs to $K$,
and the intersection of every two simplices
belonging to $K$ is a face of each of them.
The advantage of this approach is that one can
use barycentric subdivision $K'$ of $K$ obtained
by starring of $K$ at algebraic centers of its simplices $\Delta$,
and one has $\bigcup K'=\bigcup K$,
i.e. the carriers $|K|$ of $K$ and $|K'|$ of $K'$
are identical.
\par Given a geometric simplicial complex $K$ one has
a natural partition of unity on $|K|$, namely
the set of barycentric coordinates $\phi_v$, where
$v$ ranges over all vertices of $K$.
To find $\phi_v(x)$ one picks any simplex $\Delta$ of $K$
containing $x$, expresses $x$ as the linear combination
of the vertices of $\Delta$, and $\phi_v(x)$ is the coefficient
by $v$ (that means $\phi_v(x)=0$ if $v$ is not a vertex of $\Delta$).
Thus, $x=\sum\limits_{v\in V}\phi_v(x)\cdot v$, where $V$
is the set of vertices of $V$.
\par Since $|K|=|K'|$, there are two natural partitions of unity
on $|K|$ and the next result reveals the basic connection between barycentric
subdivisions and derivatives of partitions.


\begin{Prop} \label{XXX6.4} Suppose $K$ is a simplicial complex
and let, for each $v\in K^{(0)}$, $\phi_v$ be the $v$-th barycentric coordinate.
The derivative of $\{\phi_v\}_{v\in K^{(0)}}$ forms
the barycentric coordinates of the barycentric subdivision $K'$ of $K$. 

\end{Prop}

\begin{pf}  This follows from \ref{XXX4.2} and Lemmata 7-8 in \cite{M-S} (pp.306--7).
Lemma 7 can be interpreted as saying that $\{\phi_v\}_{v\in K^{(0)}}$
is the integral of $\{\phi'_\Delta\}_{v\in (K')^{(0)}}$,
and Lemma 8 gives the formula identical with that in \ref{XXX4.2}.
\end{pf}

The carrier $|K|$ of each geometric simplicial complex
can be metrized by the metric $d(x,y):=\sum\limits_{v\in V}|\phi_v(x)-\phi_v(y)|$
and the resulting metric space is denoted by $|K|_m$
(see \cite{M-S}, p.301). Since
$|K|=|K'|$, we have two metrics on the same carrier.
In traditional approaches to simplicial complexes
it is a non-trivial task to show that they are equivalent
(see \cite{M-S}, Theorem 13 on p.306).
In our approach it is a simple consequence of \ref{XXX6.4}.


\begin{Cor} \label{XXX6.5} Suppose $K$ is a simplicial complex
and $K'$ is its barycentric subdivision. The identity function $|K'|_m\to |K|_m$
is a homeomorphism. 

\end{Cor}

\begin{pf}  Proposition \ref{XXX5.4} in \cite{D$_1$}
implies that $f:X\to |K|_m$ is continuous if and only if
$\phi_v\circ f$ is continuous for each vertex $v$. 
Notice that the derivative of $\{\phi_v\circ f\}_{v\in V}$
is exactly $\{\phi_\Delta'\circ f\}_{\Delta\in K}$,
so $f:X\to |K|_m$ is continuous if and only if $f:X\to |K'|_m$
is continuous.
\end{pf}

Another easy consequence of our results is the fact that
metric simplicial complexes are absolute
neighborhood extensors of metrizable spaces which is normally
proved via Dugundji Theorem plus some non-trivial calculations
(see \cite{M-S}, Theorem 11 on p.304).
We are going to prove a stronger result.


\begin{Cor} \label{XXX6.6} Suppose $K$ is a simplicial complex,
$A$ is a closed G$_\delta$-subset of a paracompact space $X$,
and $f:A\to |K|_m$ is a continuous function. There is an open
subset $U$ of $X$ containing $A$ and a continuous extension
$g:U\to |K|_m$ of $f$.

\end{Cor}

\begin{pf}   Let $L$ be the full simplicial complex containing $K$
(that means $L$ has the same set of vertices as $K$ and
contains all possible simplices). 
Think of $f$ as a partition of unity enumerated by vertices of $L$.
Obviously, it is point-finite.
\ref{XXX0.8} says that $f$ can be extended to a point finite
partition of unity $h$ on $X$. That $h$ can be interpreted
as a continuous function from $X$ to $|L|_m$ which extends $f$.
By \ref{XXX6.3} there is a neighborhood $U$ of $A$ in $X$
and a partition of unity $g$ on $U$ extending $f|A$
such that $g$ approximates $f|U$ and its nerve equals
the nerve of $f|A$. The fact that $g$ approximates $f|U$
implies that it is point finite and can be interpreted as
a continuous function from $U$ to its nerve. That nerve is contained in $K$
(it equals the nerve of $f|A$), so one
gets an extension $g:U\to |K|_m$ of $f$.
\end{pf}

It is traditional to show that continuous functions to simplicial complexes
are homotopic if sufficiently close (see \cite{D-S} or \cite{M-S}).
Let us show that using derivatives of partitions of unity
one gets a simpler result which is reminiscent
of the well-known fact that any two continuous functions $f,g:X\to S^n$
are homotopic if $|f(x)-g(x)|<2$ for each $x\in X$.


\begin{Cor} \label{XXX6.7} Suppose $K$ is a simplicial complex
and $f,g:X\to |K|_m$ are two continuous functions which agree
on a subset $A$ of $X$. If the distance between their derivatives
$f'$ and $g'$ is less than $2$, then $f$ is homotopic to $g$
rel. $A$.

\end{Cor}

\begin{pf}  Consider the full complex $L$ containing $K$.
The identity function $id:|L|_m\to |L|_m$ may be viewed as a
partition of unity on $|L|_m$. Its derivative $(id)'$ is a partition of
unity on $|L|_m$ indexed by simplices $\Delta$ of $L$.
Also, for any continuous function $u:X\to |L|_m$,
thought of as a partition of unity on $X$,
the derivative $u'$ of $u$ equals $id'\circ u$.
Let 
$$U=\{x\in |L|_m\mid (id)'_\Delta(x)>0 \text{ and } (id)'_\Delta|f(A)\ne 0 \text{ for some }\Delta\in K\}.$$ 
Let us show that $h=(1-a)\cdot f+a\cdot g$ maps $X$ to $U$
for each $a\in [0,1]$. Since there is a retraction $r:U\to |K|_m$
(see \ref{XXX6.3}), that would complete the proof.
\par Given $x\in X$ there is $\Delta\in K$ such that
$f'_\Delta(x)>0$ and $g'_\Delta(x)>0$ (otherwise $|f'(x)-g'(x)|=2$).
If $h'_\Delta(x)=0$, then there is $s\in\Delta$ and $t\in S\setminus\Delta$
with $h_s(x)\leq h_t(x)$. However, $f_s(x)>f_t(x)$
and $g_s(x)>g_t(x)$ implies $h_s(x)=(1-a)\cdot f_s(x)+a\cdot g_s(x)>
(1-a)\cdot f_t(x)+a\cdot g_t(x)=h_t(x)$, a contradiction.
\end{pf}

We will now formalize an operation which we have already used without
mentioning it explicitely.


\begin{Prop} \label{XXX6.8} Suppose $f=\{f_s\}_{s\in S}$
is a partition of unity on an open subset $U$ of $X$
and $g=\{g_s\}_{s\in S}$
is a partition of unity on an open subset $V$ of $X$.
If $X=U\cup V$ and $\alpha:X\to [0,1]$ is a continuous
function such that $\alpha^{-1}(0,1]\subset U$
and $\alpha^{-1}[0,1)\subset V$,
then $h_s:=\alpha\cdot f_s+(1-\alpha)\cdot g_s$
defines a partition of unity on $X$ called the {\bf join} of $f$ and $g$
along $\alpha$ and denoted $f\ast_\alpha g$.

\end{Prop}

\begin{pf}  Notice that the formula for $h_s$ does not depend
on how $f_s$ is extended over $X\setminus U$
and on how $g_s$ is extended over $X\setminus V$.
The easiest choice is to extend them trivially by mapping
those complements to $0$. Applying \ref{XXX3.1}
one gets that $\{\alpha\cdot f_s\}_{s\in S}$
is a partition of $\alpha$ on $X$ and $\{(1-\alpha)\cdot g_s\}_{s\in S}$
is a partition of $1-\alpha$ on $X$ which implies that 
$\{h_s\}_{s\in S}$
is a partition of unity on $X$.
\end{pf}


\begin{Cor} \label{XXX6.9} Suppose $K$ is a simplicial complex,
$A$ is a closed G$_\delta$-subset of a paracompact space $X$,
and $f:A\to |K|_m$ is a continuous function. If there is a continuous function
$g:X\to |K|_m$ such that $g|A$ approximates $f$, then $f$ extends
continuously over $X$.

\end{Cor}

\begin{pf}  Let $S$ be the set of vertices of $K$ and let us
interpret $g$ as a partition of unity $\{g_s\}_{s\in S}$ on $X$.
Choose, using \ref{XXX1.12}, a locally finite partition of unity
$h=\{h_s\}_{s\in S}$ on $X$ such that $B_s:=cl(h_s^{-1}(0,1])
\subset g_s^{-1}(0,1]$ for each $s\in S$.
Apply \ref{XXX6.6} and find an extension $F=\{F_s\}_{s\in S}:W\to |K|_m$
of $f$ over an open neighborhood $W$ of $A$ in $X$.
Let $C_s:=F_s^{-1}(0)$ for each $s\in S$.
Notice that $A\cap C_s\cap B_s=\emptyset$ for all $s\in S$:
$x\in A\cap B_s$ implies $g_s(x)>0$ which implies $f_s(x)>0$, so $x\notin C_s$.
Since $\{C_s\cap B_s\}_{s\in S}$ is a locally finite family of closed sets in $W$,
there is an open neighborhood $U$ of $A$ in $W$
such that $U\cap C_s\cap B_s=\emptyset$ for all $s\in S$ which implies
that $h|U$ approximates $F|U$.
Pick a continuous function $\alpha:X\to [0,1]$ such that $\alpha(A)\subseteq \{1\}$
and $\alpha(X\setminus U)\subseteq \{0\}$.
Now $H:=(F|U)\ast_\alpha h=\alpha\cdot (F|U)+(1-\alpha)\cdot h$
is an extension of $f$, so it remains to show that its image is
contained in $|K|_m$, i.e. $H_s(x)>0$ for $s\in T$ implies $T\in K$.
It is certainly so for $x\in X\setminus U$ as $H_s(x)>0$ implies
$g_s(x)>0$ in that case.
If $x\in U$, then $h_s(x)>0$ implies $F_s(x)>0$ as $h|U$
approximates $F|U$, so this case holds as well.
\end{pf}

Let us show an application of \ref{XXX6.9} to the theory of absolute extensors.


\begin{Cor} \label{XXX6.10} Suppose $X$ is a metrizable space and
$K$ is a simplicial complex. If $|K|_m$ is an absolute extensor
of $X$, then it is an absolute extensor of every subset of $X$.

\end{Cor}

\begin{pf}  Case 1. Open subsets of $X$. Given an open subset $U$
of $X$ let us express it as the union $\bigcup\limits_{n=1}^\infty B_n$
of closed subsets $B_n$ of $X$ such that $B_{n}\subset int(B_{n+1})$
for all $n$. Suppose $C$ is a closed subset of $U$ and
$f:C\to |K|_m$ is a continuous function. Given an extension $f_n:B_n\to |K|_m$
of $f|B_n\cap C$, we extend $f_n$ to $f_{n+1}:B_{n+1}\to |K|_m$
so that $f_{n+1}|B_{n+1}\cap C=f|B_{n+1}\cap C$.
The direct limit of $f_n$ is an extension of $f$ over $A$.
\par Case 2: All subsets of $X$. Suppose $C$ is a closed subset of $A\subset X$ and
$f:C\to |K|_m$ is a continuous function. According to \ref{XXX6.9} it suffices to show
that an approximate of $f$ extends over $A$, so we may assume (see \ref{XXX1.12})
that $f$ is locally finite. Extend $f$ to a locally finite
partition of unity $g$ on $A$ (see \ref{XXX3.4}). 
Using \ref{XXX5.8} find a neighborhood $U$ of $A$ in $X$ and
a locally finite partition of unity $h$ on $U$ such that $h|A$ approximates $g$.
In particular
$h$ can be interpreted as a continuous function from $U$ to $|L|_m$,
where $L$ is the full simplicial complex containing $K$.
Let $D:=h^{-1}(|K|_m)$. Notice that $D$ is closed in $U$
and contains $C$ as $h|C$ approximates $f$ which implies
$h(C)\subset |K|_m$.
By Case 1, $h|D:D\to |K|_m$ extends over $U$, so $f$ extends over $A$
by \ref{XXX6.9}.
\end{pf}

Let us show that the operation of taking joins of partitions of unity
corresponds to the operation of taking joins of simplicial complexes.


\begin{Def} \label{XXX6.11} Suppose $K$ and $L$ are two abstract
simplicial complexes with sets of vertices $S_K$ and $S_L$
so that $S_K\cap S_L=\emptyset$.
The {\bf join} $K\ast L$ of $K$ and $L$ is the simplicial
complex with the set of vertices equal to $S_K\cup S_L$
so that $T\in K\ast L$ if and only if $T\cap S_K\in K$
and $T\cap S_L\in L$.

\end{Def}

Geometrically, it amounts to placing $K$ and $L$ in two
linear subspaces $E_K$ and $E_L$, respectively, of a vector space
$E$ so that $E_K\cap E_L=0$. The geometric simplices
of $K\ast L$ are obtained as convex hulls of $\sigma\cup \tau$,
where $\sigma\in K$ and $\tau\in L$.


\begin{Prop} \label{XXX6.12} Let $X$ be a topological
space and let $K,L$ be simplicial complexes
with disjoint sets of vertices $S_K$ and $S_L$, respectively.
Given a continuous function $h:X\to |K\ast L|_m$
there are a continuous function $\alpha:X\to [0,1]$
and continuous functions $f:\alpha^{-1}(0,1]\to |K|_m$,
$g:\alpha^{-1}[0,1)\to |L|_m$ such that
$h=f\ast_\alpha g$.
Conversely, given a continuous function $\alpha:X\to [0,1]$
and continuous functions $f:\alpha^{-1}(0,1]\to |K|_m$,
$g:\alpha^{-1}[0,1)\to |L|_m$,
$h=f\ast_\alpha g$ maps $X$ to $|K\ast L|_m$.

\end{Prop}

\begin{pf}  Suppose $h:X\to |K\ast L|_m$ is a continuous function.
Interpret it as a partition of unity $\{h_s\}_{s\in S}$,
$S=S_K\cup S_L$, on $X$.
Put $\alpha:=\sum\limits_{s\in S_K}h_s$ (it is continuous by \ref{XXX1.5}),
$f_s:=h_s/\alpha$ if $s\in S_K$, $f_s:=0$ if $s\in S_L$,
$g_s:=h_s/(1-\alpha)$ if $s\in S_L$, $g_s:=0$ if $s\in S_K$.
Notice that $h=f\ast_\alpha g$.
\par If $h=f\ast_\alpha g$, where 
$\alpha:X\to [0,1]$, $f:\alpha^{-1}(0,1]\to |K|_m$, and
$g:\alpha^{-1}[0,1)\to |L|_m$,
then $h_s(x)>0$ for $s\in T$ means $f_s(x)>0$ for $s\in T\cap S_K$
and $g_s(x)>0$ for $s\in T\cap S_L$, i.e. the nerve of $h$
is contained in $K\ast L$.
\end{pf}


\begin{Cor}[\cite{D$_4$}] \label{XXX6.13}   Let $X$ be a metrizable
space and let $K,L$ be simplicial complexes.
If $X=A\cup B$, $|K|_m$ is an absolute extensor of $A$,
and $|L|_m$ is an absolute extensor of $B$,
then $|K\ast L|_m$ is an absolute extensor of $X$.
 
\end{Cor}

\begin{pf} 
 Suppose $C$ is a closed subset of $X$
and $f:C\to |K*L|_m$ is a continuous function.  By \ref{XXX6.12} $f$ defines two closed,
disjoint subsets $C_K$,  $C_L$ of
$C$ and 
 continuous functions $f_K:C-C_L\to K$, $f_L:C-C_K\to L$, $\alpha :C\to
[0,1]$ such that: 
\par{1.} $\alpha ^{-1} (1)=C_K$, $\alpha ^{-1}
(0)=C_L$, 
\par{2.} $f(x)=\alpha (x)\cdot f_K(x)+(1-\alpha (x))\cdot f_L(x)$ for
all $x\in C$.
\par Since $|K|_m$ is an absolute extensor of $A-C_L$ by  \ref{XXX6.10}, 
$f_K$ extends over $(C\cup A)-C_L$. 
 Consider an approximate extension  $g_K:U_A\to K$ of $f_K$ over a
neighborhood $U_A$ of $(C\cup A)-C_L$ in $X-C_L$. 
Such an extension exists by \ref{XXX5.8} and \ref{XXX6.3}. Since
$C-C_L$ is closed in $U_A$, we may assume that $g_K$ is an actual
extension of $f_K:C-C_L\to K$ (see \ref{XXX6.9}).
Similarly, let  $g_L:U_B\to L$ be
an extension of $f_L$ over a neighborhood $U_B$ of $(C\cup B)-C_K$ 
in $X-C_K$. Notice that $X=U_A\cup U_B$. Let $\beta
:X\to [0,1]$ be an extension of $\alpha $ such that
$\beta(X-U_B)\subset \{1\}$ and $\beta(X-U_A)\subset \{0\}$. 
 Define  $f':X\to |K*L|_m$ 
as the join $g_K\ast_\beta g_L$. Notice that $f'$ is an
extension of $f$. 
\end{pf}

\begin{Rem}  V.Toni\' c \cite{T} generalized \ref{XXX6.13} to stratifiable spaces.
\cite{D$_3$} contains a generalization of \ref{XXX6.13} to hereditarily
paracompact spaces.
 
\end{Rem}

\section{Inverse limits of compact spaces. }

One of the ways of investigating compact spaces is by mapping them
to nice spaces (polyhedra, ANRs, CW complexes).
Compact spaces $X$ are often expressed as inverse limits of simpler
spaces and one of the most popular techniques is to factor continuous functions defined 
on $X$ through terms of the inverse system (see \cite{D-S} or \cite{M-S}).
The purpose of this section is to show a simple result for
partitions of unity which can be immediately applied to continuous functions
from $X$ to finite simplicial complexes. The application
to continuous functions from $X$ to ANRs and CW complexes follows from
the fact (see \cite{Hu}) that ANRs can be approximated
by simplicial CW complexes, and continuous functions to CW complexes
have a compact image contained in a finite CW complex which
is an ANR. In short, the author believes
that \ref{XXX7.2} and its proof is the basic blueprint for all the results
of similar nature.
\par
The following is a version of equicontinuity.
Indeed, \ref{XXX2.3} says that any continuous
$f:X\times Z\to (Y,d)$ so that $Z$ is compact
has the following property: for any $\epsilon>0$ and any $a\in X$
there is a neighborhood $U$ of $a$ in $X$ such that
$d(f(x,z),f(y,z))<\epsilon$ for all $x,y\in U$.
The interpretation of \ref{XXX7.1} is that, in case of continuous functions defined
on an infinite product of compact spaces, that product
can be split into two parts allowing $U$ to be the whole $X$.


\begin{Prop} \label{XXX7.1} Suppose $J$ is a partially ordered set and $\{X_j\}_{j\in J}$
is a family of compact spaces. If $f:A\to (Y,d)$
is a continuous function from a closed subset $A$ of $\prod\limits_{j\in J}X_j$
to a metric space
and $\epsilon>0$, then there is $k\in J$
so that for any pair of points $x=\{x_j\}_{j\in J}, y=\{y_j\}_{j\in J}
\in A$ the condition $x_j=y_j$ for all $j\leq k$
implies $d(f(x),f(y))<\epsilon$. 

\end{Prop}

\begin{pf}  Fix $\epsilon>0$. For $k\in J$ let $A_k$
be the set of all $x=\{x_j\}_{j\in J}\in A$
so that there is $a_k(x)=\{y_j\}_{j\in J}\in A$ with the property
that $y_j=x_j$ for $j\leq k$ but $d(f(x),f(a_k(x)))\ge\epsilon$.
Notice that $A_l\subseteq A_k$ if $k\leq l$.
If all of $A_k$ are not empty (if one of them is empty, we are done),
then there is $z=\{z_j\}_{j\in J}$ belonging to the closure
of each $A_k$.
Pick a neighborhood $U$ of $z$ in $A$ so that $d(f(x),f(y))<\epsilon$
if $x,y\in U$. We may assume that $U=A\cap \prod\limits_{j\in J}U_j$,
where $U_j$ is open in $X_j$ for each $j\in J$ and $U_j=X_j$
for all but finitely many $j\in J$.
Such
$U$ has the property that, for some $k\in J$,
$p=\{p_j\}_{j\in J}\in U$ implies $q=\{q_j\}_{j\in J}\in U$
provided $q_j=p_j$ for all $j\leq k$ and $q\in A$. Well, pick $p\in A_k\cap U$
and put $q=a_k(p)$ to arrive at a contradiction.
\end{pf}

If $X$ is the inverse limit of an inverse system
$\{X_j,p^j_i,J\}$ of topological spaces, then $p_j$ denotes the natural
projection $X\to X_j$.


\begin{Cor} \label{XXX7.2}  Suppose $(X,Y)$ is the inverse limit of an inverse system
$\{(X_j,Y_j),p^j_i,J\}$ of compact Hausdorff pairs and $\epsilon>0$.
Given a partition of unity  $f=\{f_s\}_{s\in S}$ on $Y_i$
and given an extension  $g=\{g_s\}_{s\in S}$
of $f\circ p_i$ over $X$ there is $n>i$
and an extension  $h=\{h_s\}_{s\in S}$
of $f\circ p^n_i$ over $X_n$ such that
$|h\circ p_n-g|<\epsilon$ and $|(h\circ p_n)'-g'|<\epsilon$.

\end{Cor}

\begin{pf}  Consider the subset $$Z_i=\{x=\{x_j\}_{j\in J}
\in \prod\limits_{j\in J}X_j\mid x_i\in Y_i \text{ and } x_j=p^i_j(x_i)
\text{ for all } j<i\}.$$
The projection $\pi_i:Z_i\to Y_i$ gives rise to the partition
of unity $f\circ \pi_i$ which agrees with $g$ on $X\cap Z_i$.
We can paste them together and then extend over the whole $\prod\limits_{j\in J}X_j$ using \ref{XXX0.8}.
Call the resulting partition of unity $u=\{u_s\}_{s\in S}$
and find $n\in J$, $n>i$, so that for any pair of points
$x=\{x_j\}_{j\in J}, y=\{y_j\}_{j\in J}
\in \prod\limits_{j\in J}X_j$ the condition $x_j=y_j$ for all $j\leq n$
implies $|u(x)-u(y)|<\epsilon$ and $|u'(x)-u'(y)|<\epsilon$. 
Pick points $b_j\in X_j$ for each $j\in J$ and let
 $i_n:X_n\to \prod\limits_{j\in J}X_j$ be defined
as follows: $i_n(x)=\{y_j\}_{j\in J}$,
where $y_n=x$, $y_j=p^n_j(x)$ for $j\leq n$, and $y_j=b_j$ otherwise.
Notice that $h=u\circ i_n$ satisfies the desired conditions. 
\end{pf}


\begin{Cor} \label{XXX7.3}  Let $K$ be a simplicial complex and
let $(X,Y)$ be the inverse limit of an inverse system
$\{(X_j,Y_j),p^j_i,J\}$ of compact Hausdorff pairs 
such that one of the following conditions holds:
\par{1.}
$|K|_m$ is complete,
\par{2.} $J$ is countable and each $X_j$ is compact metric.
\par\noindent
Given $\epsilon>0$, given a continuous function  $f:Y_i\to |K|_m$,
and given an extension  $g:X\to |K|_m$
of $f\circ p_i$ there is $n>i$
and an extension  $h:X_n\to |K|_m$
of $f\circ p^n_i$ such that
$|h\circ p_n-g|<\epsilon$ and $|(h\circ p_n)'-g'|<\epsilon$.

\end{Cor}

\begin{pf}  Let $L$ be the full simplicial complex containing
$K$. Let $S$ be the set of vertices of $L$.
 $l_S^1$ is the space of all functions $u:S\to R$
which are absolutely summable. All partitions
of unity $\{f_s\}_{s\in S}$ on $X$ can be viewed as continuous functions from $X$ to $l^1_S$ and
all the continuous functions to $|K|_m$ can be viewed as partitions of unity.
Consider the subset $$Z_i=\{x=\{x_j\}_{j\in J}
\in \prod\limits_{j\in J}X_j\mid x_i\in Y_i \text{ and } x_j=p^i_j(x_i)
\text{ for all } j<i\}.$$
The projection $\pi_i:Z_i\to Y_i$ gives rise to the partition
of unity $f\circ \pi_i$ which agrees with $g$ on $X\cap Z_i$.
We can paste them together and then extend over the whole $\prod\limits_{j\in J}X_j$ using \ref{XXX0.8}.
Call the resulting partition of unity $u=\{u_s\}_{s\in S}$.
$u$ maps $\prod\limits_{j\in J}X_j$ to $|L|_m$ in case 2)
and to $l^1_S$ in case 1). Indeed, in case 1) we may invoke \ref{XXX0.2}
and in case 2) we may invoke the point-finite case of \ref{XXX0.8}.
By \ref{XXX6.6} there is a retraction $r:N\to |K|_m$ from
a closed neighborhood $N$ of $|K|_m$ in $|L|_m$ (in $l^1_S$, respectively).
Let $A=u^{-1}(N)$. $int(A)$ contains $X\cup Z_i$.
By \ref{XXX7.1}, find $m\in J$, $m>i$, so that for any pair of points
$x=\{x_j\}_{j\in J}, y=\{y_j\}_{j\in J}
\in A$ the condition $x_j=y_j$ for all $j\leq m$
implies $|r\circ u(x)-r\circ u(y)|<\epsilon$ and 
$|(r\circ u)'(x)-(r\circ u)'(y)|<\epsilon$.
Let
$$B_k=\{x=\{x_j\}_{j\in J}
\in \prod\limits_{j\in J}X_j\mid x_j=p^k_j(x_k)
\text{ for all } j<k\}.$$
 We need $B_n\subset int(A)$ for some $n\ge m$.
To prove this, put $C_p:=B_p\setminus int(A)$ for $p\ge m$.
Since $\bigcap\limits_{p\ge m}B_p=X$, 
$\bigcap\limits_{p\ge m}C_p$ must be empty
and there is a finite $T\subset J$ such that
$\bigcap\limits_{p\in T}C_p=\emptyset$.
If $n$ is bigger than all elements of $T$, then $C_n=\emptyset$.
This shows $B_n\subset int(A)$.
Pick points $b_j\in X_j$ for each $j\in J$ and let
 $i_n:X_n\to \prod\limits_{j\in J}X_j$ be defined
as follows: $i_n(x)=\{y_j\}_{j\in J}$,
where $y_n=x$, $y_j=p^n_j(x)$ for $j\leq n$, and $y_j=b_j$ otherwise.
$i_n$ satisfies $i_n(X_n)\subset B_n\subset A$,
so $h=r\circ u\circ i_n$ is well-defined.
If $x\in X$, then $i_n(p_n(x))$ and $x$ have the same
coordinates up to $n$, so
$|r\circ u(i_n(p_n(x)))-r\circ u(x)|<\epsilon$ and 
$|(r\circ u)'(i_n(p_n(x)))-(r\circ u)'(x)|<\epsilon$.
Since $r(u(x))=u(x)=g(x)$ for $x\in X$,
 $h=r\circ u\circ i_n$ satisfies the desired conditions. 
\end{pf}

\section{Appendix. }

The purpose of the Appendix is to show that \ref{XXX2.3}
describes a generic way of obtaining equicontinuous
families with values in compact metric spaces. As a consequence
we get a simple proof of Ascoli Theorem.
It seems to the author that one gets a better understanding of the Ascoli Theorem
if \ref{XXX8.1} and \ref{XXX8.7} are proved first, the functorial properties
of the compact-open topology are established next,
and, finally, those properties are used to prove the result
as in \ref{XXX8.2}. By the functorial property we mean the fact
that, for $k$-spaces $X\times Z$, a function $f:Z\to Map(X,Y)$
is continuous if and only if the adjoint
function $f':X\times Z\to Y$ is continuous
(see \cite{En$_1$}, 3.4.9).


\begin{Thm} \label{XXX8.1} Let $\{f_s\}_{s\in S}$
be family of functions from a space $X$ to
a metric space $(Y,d)$. The following conditions are equivalent:
\par{a.} $\{f_s\}_{s\in S}$ is equicontinuous
and
for each $x\in X$ there is a compact subset $Y_x$
containing all values $f_s(x)$, $s$ ranging through all of $S$.
\par{b.} There is a compact Hausdorff space $Z$ and a 
continuous function $f:X\times Z\to Y$ such that the family
$\{f_z\}_{z\in Z}$ defined by $f_z(x)=f(x,z)$
contains all functions $f_s$, $s\in S$.

\end{Thm}

\begin{pf}  a)$\implies$ b).  For each $x\in X$ and for each $n>0$
let us pick a neighborhood $U(x,n)$ of $x$ in $X$
such that $d(f_s(y),f_s(x))<1/n$ for all $y\in U(x,n)$.
 Consider the set $Z$ of all
functions $g:X\to Y$ such that $d(g(y),g(x))\leq 1/n$ if $y\in U(x,n)$
and $g(x)\in Y_x$ for all $x\in X$.
Notice that $Z$ is equicontinuous by definition and
 $Z$ is a subset of $\prod\limits_{x\in X}Y_x$.
Now, give $\prod\limits_{x\in X}Y_x$ the product
topology and give $Z$ the subspace topology.
Notice that $f:X\times Z\to Y$
given by $f(x,z)=z(x)$ is continuous.
Indeed, if $(x,z)\in X\times Z$ and $n\ge 1$,
then $f^{-1}$ of the open ball around $z(x)$
of radius $1/n$ contains $U(x,2n)\times \{t\in Z\mid d(t(x),z(x))<1/(2n)\}$.
All that remains to be shown is that $Z$ is closed.
Let $g:X\to Y$ with $g\notin Z$.
Since $\prod\limits_{x\in X}Y_x$ is closed, we may assume $g\in \prod\limits_{x\in X}Y_x\setminus Z$.
Suppose $d(g(b),g(a))>1/n$
for some $b\in U(a,n)$ and put $\epsilon=(d(g(b),g(a))-1/n)/3$.
Let $V$ be the set of all functions $h:X\to Y$
such that $d(h(b),g(b))<\epsilon$ and $d(h(a),g(a))<\epsilon$.
$V\cap\prod\limits_{x\in X}Y_x$ is an open subset of $\prod\limits_{x\in X}Y_x$
and is contained in $\prod\limits_{x\in X}Y_x\setminus Z$.
Indeed,
$d(h(b),h(a))\ge d(g(b),g(a))-d(h(a),g(a))-d(h(b),g(b))>1/n$.
\par b)$\implies$ a). This follows from \ref{XXX2.3} and the fact that
$Y_x=f(\{x\}\times Z)$ is compact.
\end{pf}

Recall that if $\cal F\subseteq Map(X,Y)$
is a subset of functions from $X$ to $Y$,
then we have a natural function called {\bf the evaluation map}
$eval:X\times \cal F\to Y$ defined by $eval(x,f)=f(x)$.


\begin{Thm}[Ascoli \cite{En$_1$}, 3.4.20] \label{XXX8.2} 
Let $X$ be a $k$-space.
Suppose $Y$ is a metric space
and $\cal F\subseteq Map(X,Y)$ is a subspace
of the space of continuous functions from $X$ to $Y$
considered with the compact-open topology.
If $\cal F$ is equicontinuous and $eval(\{x\}\times\cal F)$
is contained in a compact subset of $Y$ for each $x\in X$, then 
the closure of $\cal F$ in
$Map(X,Y)$ is compact.

\end{Thm}

\begin{pf}  By \ref{XXX8.1}, pick a continuous function $f:X\times Z\to Y$
such that $Z$ is compact Hausdorff
and $\cal F$ is contained in $\{f_z\}_{z\in Z}$.
Since $X$ is a $k$-space, $X\times Z$ is a $k$-space and the induced
function $g:Z\to Map(X,Y)$, $g(z)(x)=f(x,z)$, is continuous
(see \cite{En$_1$}, 3.3.27).
Notice that $g(Z)$ contains $\cal F$. 
\end{pf}


\begin{Rem} \label{XXX8.3} 
In \cite{D$_2$} (see Theorem 4.17) the author
stated an Ascoli Type Theorem involving the so-called
covariant topology on function spaces introduced there.
It dealt with $k$-spaces as in \ref{XXX8.2}.
It is clear now that $X$ does not have to be a $k$-space at all
(the function $g$ in the above proof is always continuous
if $Map(X,Y)$ is given the covariant topology)
which indicates that the covariant topology makes sense.

\end{Rem}

What should be the meaning of the concept
of equicontinuity of $\cal F\subseteq Map(X,Y)$
for arbitrary, not necessarily metric, $Y$?
The author believes that the answer ought to be
as follows.


\begin{Def}[Heuristic Definition] \label{XXX8.4} $\cal F$ is {\bf equicontinuous} if
there is a compact space $Y'$ containing $Y$
and there is an extension $f:X\times Z\to Y'$
of the evaluation function $eval:X\times \cal F\to Y'$
such that $f$ is continuous, $Z$ is compact,
and $Z$ contains $\cal F$ as a subset. 

\end{Def}

We will show that \ref{XXX8.4} makes sense in the case of completely regular $Y$.
\par First, recall the definition of equicontinuity
from \cite{En$_1$}, 3.4.17-20.


\begin{Def} \label{XXX8.5} Let $X$ and $Y$ be topological spaces.
A family $\{f_s:X\to Y\}_{s\in S}$ is {\bf equicontinuous} if
for each $x\in X$, each $y\in Y$, and each neighborhood $V$ of $y$
in $Y$ there exist neighborhoods $U$ of $x$ in $X$ and $W$ of $y$
in $Y$ such that, for every $s\in S$, $f_s(x)\in W$ implies
$f_s(U)\subseteq V$.

\end{Def}

Beware of the fact that \ref{XXX1.7} deals with functions to a space with a specified
metric and \ref{XXX8.5} deals with functions to topological spaces.
It is easy to check that if $\{f_s:X\to (Y,d)\}_{s\in S}$ is equicontinuous
in the sense of \ref{XXX1.7}, then $\{f_s:X\to Y\}_{s\in S}$ is equicontinuous
in the sense of \ref{XXX8.5}. The converse may not be true:
Consider $X=(0,1]=Y$ and $f_n(x):=x/n$ for $n\ge 1$.
$\{f_n:X\to (Y,d)\}_{n\ge 1}$ is equicontinuous
in the sense of \ref{XXX1.7} if $d$ is the standard metric ($d(a,b)=|a-b|$)
but is not equicontinuous
in the sense of \ref{XXX1.7} if $d(a,b):=|1/a-1/b|$.
However, \ref{XXX8.1} and \ref{XXX8.6}-\ref{XXX8.7} show that \ref{XXX1.7} and \ref{XXX8.5}
are equivalent for families of functions $\{f_s:X\to Y\}_{s\in S}$
such that for each $x\in X$ there is a compact subset $Y_x$ of $Y$
containing all values $f_s(x)$, $s$ ranging through all of $S$.
That is a very important class of functions in view of applications
via the Ascoli Theorem.


\begin{Lem} \label{XXX8.6} Suppose $f:X\times Z\to Y'$ is continuous
and $Y\subseteq Y'$. If $Z$ is compact, $Y'$ is regular,
and $S:=\{s\in Z\mid f(X\times\{s\})\subseteq Y\}$,
then the induced family of functions
$\{f_s:X\to Y\}_{s\in S}$, $f_s(x)=f(x,s)$, is
equicontinuous.

\end{Lem}

\begin{pf}  Suppose $x\in X$, $y\in Y$, and $V$ is a neighborhood of $y$ in $Y$.
Pick an open set $V'$ in $Y'$ satisfying $V=Y\cap Y'$.
For every pair $(U,W)$ such that $U$ is a neighborhood of $x$ in $X$
and $W$ is a neighborhood of $y$ in $Y$ define
$A(U,W):=\{s\in S\mid f(x,s)\in W \text{ and }f(U\times\{s\})\setminus V\ne\emptyset\}$.
It suffices to show $A(U,W)=\emptyset$ for some $U$ and some $W$.
Suppose, on the contrary, that none of those sets is empty.
Since $A(U',W')\subseteq A(U,W)$ if $U'\subseteq U$
and $W'\subseteq W$,
there is $z_0\in Z$ belonging to closures of all $A(U,W)$'s.
\par
 If $y\ne f(x,z_0)$, then we pick a neighborhood $W$ of $y$ in $Y'$
whose closure $cl(W)$ misses $f(x,z_0)$.
Thus, $f(x,z_0)\in Y'\setminus cl(W)$
and there is a neighborhood $U\times U'$ of $(x,z_0)$
in $X\times Z$ such that $f(U\times U')\subseteq V'\setminus cl(W)$.
Therefore there is $s\in U'\cap A(U,W\cap Y)$.
However, $s\in U'$ implies $f(U\times\{s\})\subseteq V'\setminus cl(W)$,
and $s\in A(U,W\cap Y)$ implies $f(x,s)\in W\cap Y$, a contradiction.
\par Thus, $f(x,z_0)=y\in V'$ and there is a neighborhood $U\times U'$ of $(x,z_0)$
in $X\times Z$ such that $f(U\times U')\subseteq V'$.
Therefore there is $s\in U'\cap A(U,Y)$.
However, $s\in U'$ implies $f(U\times\{s\})\subseteq V'\cap Y=V$,
and $s\in A(U,Y)$ implies $f(U\times\{s\})\setminus V\ne\emptyset$, a contradiction.
\end{pf}


\begin{Thm} \label{XXX8.7} Let $\{f_s\}_{s\in S}$
be an equicontinuous family of functions from a space $X$ to
a regular space $Y$. If 
for each $x\in X$ there is a compact subset $Y_x$ of $Y$
containing all values $f_s(x)$, $s$ ranging through all of $S$,
then there is a compact Hausdorff space $Z$ and a 
continuous function $f:X\times Z\to Y$ such that the family
$\{f_z\}_{z\in Z}$ defined by $f_z(x)=f(x,z)$
contains all functions $f_s$, $s\in S$.

\end{Thm}

\begin{pf}  For each triple $(x,y,V)$ such that $(x,y)\in X\times Y$
and $V$ is a neighborhood of $y$ in $Y$
pick a neighborhood $V'$ of $y$ in $V$ such that $cl(V')\subseteq V$,
and pick neighborhoods $U(x,y,V)$ of $x$ in $X$
and $W(x,y,V)$ of $y$ in $Y$ such that
$f_s(x)\in W(x,y,V)$ implies $f_s(U(x,y,V))\subseteq V'$.
Consider the set $Z$ of all functions $g:X\to Y$ such that
$g(x)\in W(x,y,V)$ implies $g(U(x,y,V))\subseteq cl(V')$
and $g(x)\in Y_x$ for each $x$.
 $Z$ is equicontinuous by definition.
Notice that $Z$ is a subset of $\prod\limits_{x\in X}Y_x$.
Now, give $\prod\limits_{x\in X}Y_x$ the product
topology and give $Z$ the subspace topology.
Notice that $f:X\times Z\to Y$
given by $f(x,z)=z(x)$ is continuous.
Indeed, if $(x,z)\in X\times Z$ and 
$V$ is open in $Y$, $f(x,z)\in V$,
then $f^{-1}(V)$ contains $U(x,f(x,z),V)\times \{h\in Z\mid h(x)\in W(x,f(x,z),V)\}$.
All that remains to be shown is that $Z$ is closed.
Let $g:X\to Y$ with $g\notin Z$.
Since $\prod\limits_{x\in X}Y_x$ is closed, we may assume $g\in \prod\limits_{x\in X}Y_x\setminus Z$.
Suppose $g(x')\in Y\setminus cl(V')$ and $g(x)\in W(x,y,V)$
for some $x'\in U(x,y,V)$.
Consider $\{h\in \prod\limits_{x\in X}Y_x\mid h(x')\in Y\setminus cl(V')\text{ and }h(x)\in W(x,y,V)\}$.
This is an open set missing $Z$ and containing $g$.
\end{pf}

Notice that following the proof of \ref{XXX8.2} one can give a proof
of the part of 3.4.20 in
\cite{En$_1$} which deals with proving that the closure of certain subspaces of $Map(X,Y)$
is compact. To prove that theorem completely one only
needs to use functorial properties of the compact-open
topology.
\par
The next result follows easily from \ref{XXX8.6} and \ref{XXX8.7}.


\begin{Cor} \label{XXX8.8} Let $\{f_s\}_{s\in S}$
be a family of functions from a space $X$ to
a completely regular space $Y$.
The following conditions are equivalent:
\par{1.} $\{f_s\}_{s\in S}$, $f_s(x)=f(x,s)$, is
equicontinuous. 
\par{2.} For each compact Hausdorff space $Y'$
containing $Y$ there is a compact Hausdorff space
$Z$ and a continuous function $f:X\times Z\to Y'$
such that
$\{f_z\}_{z\in Z}$ defined by $f_z(x)=f(x,z)$
contains all functions $f_s$, $s\in S$.
\par{3.} There is a compact Hausdorff space $Y'$
containing $Y$, there is a compact Hausdorff space
$Z$, and there is a continuous function $f:X\times Z\to Y'$
such that
$\{f_z\}_{z\in Z}$ defined by $f_z(x)=f(x,z)$
contains all functions $f_s$, $s\in S$.

\end{Cor}

Finally, let us explain the concept of strong equicontinuity.


\begin{Prop} \label{XXX8.9} Suppose $\{f_s:X\to [0,\infty)\}_{s\in S}$
is a family of functions on a topological space $X$.
Let $\omega(S):=S\cup\{\infty\}$ be the one-point compactification of $S$
considered with the discrete topology.
$\{f_s\}_{s\in S}$ is strongly equicontinuous if and only if
the function $f:X\times \omega(S)\to [0,\infty)$
is continuous, where $f(x,s)=f_s(x)$ for $s\in S$
and $f(x,\infty)=0$ for all $x\in X$.

\end{Prop}

\begin{pf}  $f$ is continuous at $(x,\infty)$ if and only if
for each $\epsilon>0$
there is a neighborhood $U$ of $x$ in $X$ and a neighborhood $V$ of
$\infty$ in $\omega(S)$
such that $f(U\times V)\subset [0,\epsilon)$.
That means $T:=S\setminus V$ is finite and $f_s(x)<\epsilon$ for all
$s\in S\setminus T$.
$f$ is continuous at $(x,s)$ if and only if $f_s$ is continuous at $x$.
\end{pf}

\begin{Rem}  In \cite{Y$_2$} K. Yamazaki
 proved that every pointwise bounded equicontinuous collection of 
real-valued (or Frechet-spaces- valued) functions on $A$ can be extended 
to a pointwise bounded equicontinuous collection of functions on $X$ if and only if $A$ is $P$-embedded in $X$.
 
\end{Rem}

\medskip
\medskip
\medskip
\medskip
Jerzy Dydak,
Math Dept, University of Tennessee, Knoxville, TN 37996-1300, USA,
E-mail addresses: dydak at math.utk.edu


\begin{thebibliography}{999}

\bibitem{Bor}
K.Borsuk,
{\em Theory of retracts},
Polish Scientific Publishers (1967), Warszawa.

\bibitem{Di}
Tammo tom Dieck,
{\em Partitions of unity in homotopy theory},
Compositio Mathematica 23 (1971),159--167.

\bibitem{Do}
A.Dold,
{\em Partitions of unity in the theory of fibrations},
Ann. of Math. 78 (1963),223--255.

\bibitem{D$_1$}
J.Dydak,
{\em Extension theory:  The interface between set-theoretic and algebraic topology},
Topology and its Appl. 20 (1996),1--34.

\bibitem{D$_2$}
J.Dydak,
{\em Covariant and contravariant points of view in topology with applications to function spaces},
Topology and its Appl. 94 (1999),87--125.

\bibitem{D$_3$}
J.Dydak,
{\em Extension dimension of paracompact spaces},
 to appear in Topology and its Applications.

\bibitem{D$_4$}
J.Dydak,
{\em Cohomological dimension and metrizable spaces II},
Trans.Amer.Math.Soc. 348 (1996),1647--1661.

\bibitem{D-F}
J.Dydak and N.Feldman,
{\em Major theorems on compactness: A unified approach},
American Mathematical Monthly 99 (1992),220--227.

\bibitem{D-S}
J.Dydak and J.Segal,
{\em Shape theory: An introduction},
Springer Verlag (1978).

\bibitem{En$_1$}
R.Engelking,
{\em General Topology},
 (1989), Berlin.

\bibitem{En$_2$}
R.Engelking,
{\em Theory of Dimensions Finite and Infinite},
Heldermann Verlag (1995).

\bibitem{Ho}
T.Hoshina,
{\em Extensions of mappings II},
in book in Topics in General Topology, Elsevier Science Publishers, (1989),41--80.

\bibitem{Hu}
S.T.Hu,
{\em Theory of retracts},
Wayne State University Press (1965).

\bibitem{J}
I.M.James,
{\em General Topology and Homotopy Theory},
Springer Verlag (1984).

\bibitem{L-W}
A.T.Lundell and S.Weingram,
{\em The topology of CW complexes},
Van Nostrand Reinhold Co. (1969), New York.

\bibitem{M}
E.Michael,
{\em Another note on paracompact spaces},
Proceedings of the American Math.Soc. 20 (1957),41--80.

\bibitem{M-S}
S.Marde\` si\' c and J.Segal,
{\em Shape theory},
North-Holland Publ.Co. (1982), Amsterdam.

\bibitem{Mu}
James R.Munkres,
{\em Topology},
Prentice Hall (2000), 2nd edition.

\bibitem{T}
V.Toni\' c,
{\em Extension theory analogue of Menger-Urysohn Addition Theorem for stratifiable spaces},
Abstracts of Geometric Topology II Conference held in Dubrovnik (Sept 29- Oct 5, 2002),p.55.

\bibitem{Y$_1$}
K.Yamazaki,
{\em Controlling extensions of functions and C-embedding},
Topology Proceedings, to appear.

\bibitem{Y$_2$}
K.Yamazaki,
{\em Extensions of pointwise bounded equicontinuous collections of functions},
preprint.
\end{thebibliography}
\end{document}